\documentclass[eqno,11pt]{article}

\usepackage[usenames, dvipsnames]{xcolor}
\usepackage{tikz}

\usepackage{amssymb,hyperref,amsthm}
\usepackage{euscript}
\usepackage{flafter}
\usepackage{pstricks}
\usepackage[latin1]{inputenc}
\usepackage{amsmath}
\usepackage{amsfonts}
\usepackage{pstricks-add}
\usepackage{dsfont}
\usepackage{bm}
\usepackage{variations}
\usepackage{graphicx,color,dsfont,epsfig,caption,subcaption,wrapfig}
\usepackage{xcolor}
\usepackage{tikz}

\hypersetup{
    linktoc=page,
   linkcolor=red,          % color of internal links
  citecolor=blue,        % color of links to bibliography
    filecolor=blue,      % color of file links
   urlcolor=cyan,
    colorlinks=true           % color of external links
}

\usepackage[refpage]{nomencl}

\usepackage{makeidx}
\makeindex

\usepackage{mathrsfs} %%% Nice caligraphic fonts
\evensidemargin\oddsidemargin
\newcommand{\eqnsection}{
\renewcommand{\theequation}{\thesection.\arabic{equation}}
   \makeatletter
   \csname  @addtoreset\endcsname{equation}{section}
   \makeatother}
\eqnsection
\def\r{{\mathbb R}}
\def\q{{\mathbb Q}}
\def\e{{\mathbb E}}
\def\p{{\mathbb P}}
\def\P{{\bf P}}
\def\E{{\bf E}}

\def\N{{\mathbb N}}
\def\T{{\mathbb T}}
\def\L{{\mathscr  L}}
 \def\ee{\mathrm{e}}
\def\mm{{\tt m\, }} 
\newtheorem{theorem}{Theorem}[section]

\newtheorem{Lemma}[theorem]{Lemma}
\newtheorem{Proposition}[theorem]{Proposition}
\newtheorem{Remark}[theorem]{Remark}

\newtheorem{Claim}[theorem]{Claim}

\begin{document}

\vglue20pt

\begin{center}
{\Large\bf The minimum of a branching random walk \\ outside the boundary case}
\end{center}
 
\bigskip
\bigskip

\medskip

\centerline{Julien Barral, Yueyun Hu   and Thomas Madaule\footnote{J.B.  and Y.H.:   D\'epartement de Math\'ematiques, Universit\'e Paris XIII, Sorbonne Paris Cité, 99 avenue J-B Cl\'ement,  F-93430 Villetaneuse.  barral@math.univ-paris13.fr  and    yueyun@math.univ-paris13.fr 
\\ $\phantom{aob}$  T.M.:   Institut de Mathématiques de Toulouse, Université Paul Sabatier,  118, route de Narbonne,  F-31062 Toulouse Cedex 9.  thomas.madaule@math.univ-toulouse.fr 
\\ $\phantom{aob}$ Research partially  supported by ANR  (MEMEMOII)  2010 BLAN 0125}}

\medskip

\centerline{\it Universit\'e Paris XIII, \,  Universit\'e Paris XIII  \, \& Université Paul Sabatier}

\medskip

\bigskip
\bigskip

{\leftskip=2truecm \rightskip=2truecm \baselineskip=15pt \small

\noindent{\slshape\bfseries Summary.} 
 %This paper deals with a case left open in the studies of the asymptotic behavior of the minimum of a real-valued branching random walk.     
This paper is a complement to the studies on the minimum of a real-valued branching random walk.    In the boundary case  (\cite{BK05}),  A\"{i}dékon  in a seminal paper (\cite{Aidekon13}) obtained  the convergence in law of the minimum after a suitable renormalization.  We study   here   the   situation  when the log-generating function of the branching random walk explodes at some positive point and it cannot be reduced to the boundary case. In the associated  thermodynamics framework this corresponds to a first order phase transition, while the boundary case corresponds to a second order phase transition. 
%Assuming that the step distribution seen from the Peyri\`ere measure is absolutely continuous near $-\infty$ and belongs to the domain of attraction of a stable law, we prove a tightness result for the minimum, by showing that  all  particles near to the minimum are such that their  ancestors have exactly one huge drop.   If furthermore,   we   assume  the convergence in law of the point process associated to the  huge drops at  the same generation, then the renormalized minimum converges   to some    distribution convoluted by the non-degenerated limit of the additive martingale.  
\bigskip

 \noindent{\slshape\bfseries Keywords.}   Branching random walk, minimal position,   phase transition. 

\bigskip

 \noindent{\slshape\bfseries 2010 Mathematics Subject  Classification. 60J80, 60F05. }  

} %%%%%% End of narrower

\bigskip

\section{Introduction} %$\ $

%The branching random walk  is a classical model on the evolution of some population with spatial movements. 

Consider a branching random walk on the real line $\r$. Initially, a particle sits at the origin. Its
children form the first generation; their displacements from the origin correspond
to a point process $\L$ on the line. These children have children of their own (who form
the second generation), and behave,  relative  to their respective positions, like
independent copies of $\L$,  and so on. Denote by $\p$ the probability distribution on the space $\Omega$ of marked trees associated with this branching random walk, and $\mathbb E$ the expectation with respect to $\mathbb P$. 
%\medskip

The genealogy of all particles forms a Galton-Watson tree $\T$ whose root is denoted by $\varnothing$.  Denote by    $\{u : |u|=n\}$   the set of particles at generation  $n\in {\mathbb N}$  and by   $V(u) \in \r$ the position of $u$.    Notice that  $\sum_{|u | =1} \delta_{\{ V(u)\}}= \L$.  Let  $\phi$ be the log-generating function of $\L$:    $$ \phi(\beta):= \log \mathbb E\Big[ \sum_{ \vert u \vert =1} \ee^{- \beta V(u)} \Big]=  \log \mathbb E \left[ \int_\r \ee^{- \beta x}  \L(dx)\right] \in  ( - \infty, \infty], \qquad  \beta \in \r. $$

We assume that  $\T$ is  supercritical and define $M_n:= \min_{ \vert u \vert = n} V(u)$  the minimum  of the branching random walk in the $n$th generation  (with convention: $\inf\emptyset  \equiv \infty$).    Hammersley~\cite{Hammersley74}, Kingman~\cite{Kingman75} and Biggins~\cite{Biggins76}  have established the law of large numbers for $M_n$  under a fairly general setting: if $\mathrm{dom}(\phi)\cap\mathbb R_+^*\neq\emptyset$ then  upon on the survival of the system,  $\lim_{n\to\infty}{M_n \over n} = c$,  where $c=-\inf\{\phi(\beta)/\beta: \beta>0\}$.   Hammersley \cite{Hammersley74} raised  the problem of finding the asymptotic behavior of $M_n-cn$. Several recent attempts led to significative contributions (see \cite{AddarioReed}, \cite{BramsonZeitouni09}, \cite{HuShi09} and the references therein), until the  sharp answer was given by A\"\i d\'ekon in~\cite{Aidekon13} in the  ``boundary case'' (in the senses of \cite{BK05}, see below). 

Due to the interplay between branching random walk theory and some random energy models in statistical physics, we find useful to describe the above-mentioned fine results on $M_n$  as being obtained under a second order phase transition.  Indeed, suppose that $\mathrm{dom}(\phi)\cap\mathbb R_+^*\neq\emptyset$.  Either $c=\lim_{\beta \to\infty} -\phi(\beta)/\beta$ or $- \inf\{\phi(\beta)/\beta: \beta>0\}$ is reached at a unique $\beta_c>0$. In the latter case,  with $c$ is associated a phase transition phenomenon:  define the convex functions
$$
F_n(\beta)=\frac{1}{n}\log \sum_{|u|=n}\ee^{-\beta V(u)},\quad n\ge 1,\ \beta>0.
$$
In the random energy model introduced by  Derrida and Spohn in \cite{DerSpo} (in which $\mathbb T$ is a regular tree and the increments of the branching random walks are i.i.d. and Gaussian), these functions are the partition functions of  the directed polymers  on the disordered  tree  $\mathbb T $. They converge almost surely pointwise on  $\mathbb R_+$  to the free energy in infinite volume 
$$
F(\beta)= \mathbf{1}_{[0,\beta_c]}(\beta)\phi(\beta)+\mathbf{1}_{(\beta_c,\infty)}(\beta)\beta_c^{-1}\phi(\beta_c) \beta,\quad \beta>0, $$
(see  \cite{CK,Biggins92,Mol,OW,AB}). To slightly simplify the discussion, suppose that $\phi'(\beta_c-)$ exists (this is the case for instance when the branching number $\sum_{|u|=1}1$ has a finite expectation). When $\beta_c^{-1}\phi(\beta_c)=\phi'(\beta_c-)$, $F$ is twice differentiable everywhere except at $\beta_c$, where it is only once differentiable; in the thermodynamical setting this corresponds to a second order phase transition at temperature $\beta_c^{-1}$. When $\beta_c^{-1}\phi(\beta_c)>\phi'(\beta_c-)$, $F$ is differentiable everywhere except at $\beta_c$, and we face a first order phase transition at $\beta_c^{-1}$.  

\medskip

By using the linear transform $(V(u), u \in \T) \to  (\beta_c V(u)+\phi(\beta_c) |u|, u \in \T)$ one reduces the two previous situations  to the case where $\beta_c=1$ and 
\begin{equation}\label{Hyp1}
\phi(1)=0
\end{equation}
(see Figures~\ref{Fig1} and~\ref{Fig2}).
 \begin{figure}[ht]
\begin{subfigure}[b]{0.45\textwidth}
\centering
\begin{tikzpicture}[xscale=0.6,yscale=0.6]
\draw [->] (0,-1) -- (0,6) node [above] {$\phi(\beta)$};
\draw [->] (-2,0) -- (5,0) node [right] {$\beta$};
\draw [thick,domain=0:2] plot (\x, {3*0.25*(\x/2-4)*(\x/2-1)});
\draw [dashed,domain=-1:0] plot (\x, {3*0.25*(\x/2-4)*(\x/2-1)});
\draw [thick] (2,7) -- (5,7) node [right] {$+\infty$};
\draw [fill] (0,0) circle [radius=0.04] node [below left] {$0$};
\draw [fill] (2,0) circle [radius=0.04] node [above] {$1$};
\end{tikzpicture}
%\caption{First order phase transition: \\
%log-generating function $\phi$.}\label{fig3a}
\end{subfigure}
\begin{subfigure}[b]{0.45\textwidth}
\centering
\begin{tikzpicture}[xscale=0.6,yscale=0.6]
\draw [->] (0,-1) -- (0,6) node [above] {$F(\beta)$};
\draw [->] (-2,0) -- (5,0) node [right] {$\beta$};
\draw [thick,domain=0:2] plot (\x, {3*0.25*(\x/2-4)*(\x/2-1)});
\draw [dashed,domain=-1:0] plot (\x, {3*0.25*(\x/2-4)*(\x/2-1)});
\draw [thick] (2,0) -- (5,0); 
\draw [fill] (0,0) circle [radius=0.04] node [below left] {$0$};
\draw [fill] (2,0) circle [radius=0.04] node [above] {$1$};
\end{tikzpicture}
%\caption{First order phase transition: \\
%Free energy function $F$.}\label{fig3a}
\end{subfigure}
\caption{First order phase transition}\label{Fig1}
\end{figure}
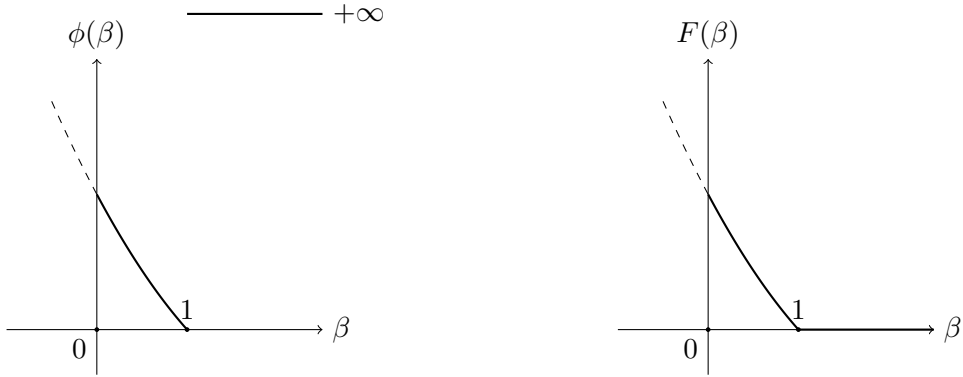

 \begin{figure}[ht]
\begin{subfigure}[b]{0.45\textwidth}
\centering
\begin{tikzpicture}[xscale=0.6,yscale=0.6]
\draw [->] (0,-1) -- (0,6) node [above] {$\phi(\beta)$};
\draw [->] (-2,0) -- (5,0) node [right] {$\beta$};
\draw [thick,domain=0:2] plot (\x, {0.75*(\x-2)^2});
\draw [dashed,domain=-0.5:0] plot (\x, {0.75*(\x-2)^2});
\draw [dashed,domain=2:3] plot (\x, {0.75*(\x-2)^2});
\draw [fill] (0,0) circle [radius=0.04] node [below left] {$0$};
\draw [fill] (2,0) circle [radius=0.04] node [above] {$1$};
\end{tikzpicture}
%\caption{First order phase transition: \\
%log-generating function $\phi$.}\label{fig3a}
\end{subfigure}
\begin{subfigure}[b]{0.45\textwidth}
\centering
\begin{tikzpicture}[xscale=0.6,yscale=0.6]
\draw [->] (0,-1) -- (0,6) node [above] {$F(\beta)$};
\draw [->] (-2,0) -- (5,0) node [right] {$\beta$};
\draw [thick,domain=0:2] plot (\x, {0.75*(\x-2)^2});
\draw [dashed,domain=-0.5:0] plot (\x, {0.75*(\x-2)^2});
\draw [thick] (2,0) -- (5,0); 
\draw [fill] (0,0) circle [radius=0.04] node [below left] {$0$};
\draw [fill] (2,0) circle [radius=0.04] node [above] {$1$};
\end{tikzpicture}
%\caption{First order phase transition: \\
%Free energy function $F$.}\label{fig3a}
\end{subfigure}
\caption{Second order phase transition (boundary case)}\label{Fig2}
\end{figure}
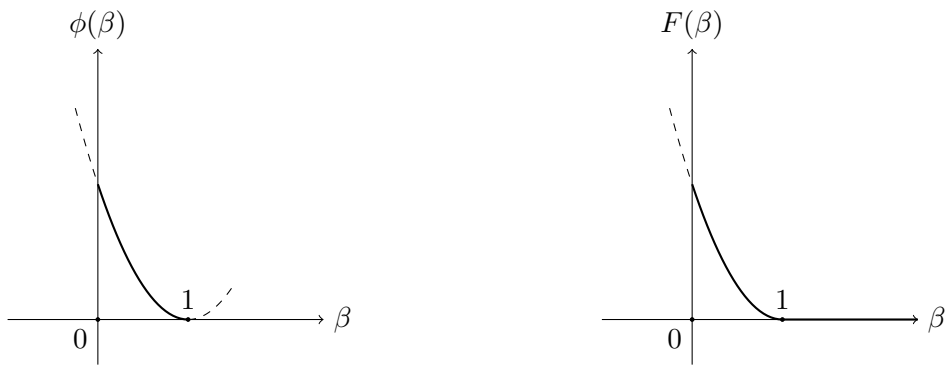

We assume \eqref{Hyp1} throughout this paper.  In case of a second order phase transition we have $\phi'(1-)=0$, i.e. $\mathbb E\big[ \int_\r  x \, \ee^{-   x} \L(dx)\big]=0$, and this last property corresponds to the  ``boundary case'' (a terminology introduced in \cite{BK05}) or the ``critical case'' in the study of the additive martingale 
$$
W_n:=\sum_{|u|=n} \ee^{-V(u)}, \qquad n\ge 1, 
$$  
while in case of a first  order phase transition we have $\phi'(1-)<0$, which more generally rewrites $\mathbb E\big[ \int_\r  x \, \ee^{-   x} \L(dx)\big] >0$, and corresponds to the so-called  ``subcritical case''. Also, since $\beta_c=1$, due to the convexity of $\phi$ we necessarily have $\phi(\beta)=+\infty$ for all $\beta>1$ in the subcritical case. Notice that in both critical and subcritical cases, the limiting velocity $c=0$. 

\medskip

When a second order phase transition occurs (namely the boundary case),  the almost sure limit of  $W_n$ vanishes (Biggins \cite{biggins-mart-cvg}, Lyons \cite{Lyons97}). Under some integrability conditions  and in the case that $\L$ is not a.s. supported on a deterministic   lattice,  it is known that the branching random walk exhibits some highly non-trivial universalities:  In the seminal paper  \cite{Aidekon13}, A{\"{\i}}d{\'e}kon proved (see also  \cite{BramsonDingZeitouni14+} for an alternative approach) the convergence in law for $M_n - {3\over 2} \log n$ as $n \to \infty$ towards a convoluted Gumbel distribution; specifically there exists a constant $c'>0$ depending on the distribution of $\L$ such that 
\begin{equation}\label{loidumin}
\lim_{n\to\infty} \mathbb P(M_n\ge \frac{3}{2}\log n +x)=\mathbb E\left( \exp(- c' \ee^x D_\infty)\right),  \quad  \forall x\in\mathbb R,
\end{equation}

\noindent where $D_\infty:= \lim_{n \to\infty} \sum_{|u|=n} V(u)\,\ee^{- V(u)}$, non-trivial and nonnegative, is the limit of the so-called derivative martingale (\cite{BigginsKyprianou04, Aidekon13, Chen2014+}). 
%We refer to Chen \cite{Chen2014+} for an optimal condition on the non-triviality of $D_\infty$ in the boundary case.  
This behavior is analogous to that observed in the branching Brownian motion (see \cite{Bramson78}).  
%  The intrinsic constant ${3\over2}$ can be explained as follows: In any   case, a certain random walk $(S_n)$   is associated with the branching random walk; the distribution of  $(S_n)$' increments will be given  in \eqref{defS}. In the boundary case,   the random walk $S_n$ has expectation $n \mathbb E\big[ \int_\r  x \, \ee^{-   x} \L(dx)\big] $, which is equal to 0 since $\phi'(1-)=0$, and it is assumed to have a finite variance, so that it makes an excursion of length $n$  with probability  $n^{-{3\over2}+o(1)}$; naturally different constants  from $\frac{3}{2}$ may appear under   weaker     moment finiteness conditions (see Remark 2 in \cite{Aidekon13}).  
%  
It is worth mentioning that  A{\"{\i}}d{\'e}kon's result \eqref{loidumin} is a key  point in understanding the asymptotic behaviors of the Gibbs measures  $\mu_{\beta,n}$  whichs assigns   to each bond $v$ of generation $n$ the mass $\mu_{\beta,n}(u)=\ee^{-\beta V(u)-n F_n(\beta)}$.    Based on  \cite{Aidekon13},  Madaule \cite{Madaule11} showed  that $n^{\frac{3}{2}\beta}\sum_{|u|=n} \ee^{-\beta V(u)}$  converges in law, see also   Webb \cite{Webb11} in the Gaussian case on a regular tree.  In the case where $\mathbb T$ is regular, say $s$-adic,   Barral, Rhodes and Vargas \cite{barhovar12} showed,  thanks to \cite{Madaule11} and the theory of invariant distributions by random weighted means (also called fixed points of the smoothing transformation theory) \cite{mandel,KP,Biggins76,DL,BK97,liu98,AlBiMe}, that   for each $\beta>1$,    $\mu_{\beta,n}$ converges in law to 
a  random discrete measure $\mu_\beta$ defined as follows:  Let $\mu$ be  the critical Mandelbrot measure on $\{0,\ldots,s-1\}^{\N_+}$ associated with the branching random walk,  that is the measure which assigns mass $\ee^{-V(u)}  D_\infty(u)$ to bond $u$, where $D_\infty(u)$ is the copy of $D_\infty$ built with the branching random walk rooted at $u$; let $N^{(\beta)}_{\mu}$ be a positive Borel random measure on $\{0,\ldots,s-1\}^{\N_+}\times \mathbb R^*_+$ whose law conditionally on $\mu$ is that of a Poisson point measure with intensity $\frac{ \mu(dx)dz}{z^{1+1/\beta}}$; then define the random measures $ \nu_{\beta} (A)= \int_A\int_{\mathbb R_+^*} z N^{(\beta)}_{\mu}(dx,dz)$ and $\mu_\beta=\nu_\beta/\|\nu_\beta\|$. All these results provide a sharp description of the asymptotic behavior of the associated directed polymer at temperatures lower than the critical freezing temperature $\beta_c=1$. In particular, they describe in which way the lower is the temperature, the more the main part of the energy concentrates on a small number of atoms.   

Let us also mention that  $M_n$ plays a role in the study of the modulus of continuity of the  0-dimensional measure $\mu$ (\cite{BKNSW}).

\medskip

In this paper we seek for the asymptotic behaviors of $M_n$ in the situation when a first order phase transition occurs, and which can not be reduced to the boundary case.   We show a convergence similar to \eqref{loidumin}  with some norming sequence depending on the law of $\L$ instead of the universal $({3\over2}\log n)$ recentering in the boundary case,  and with $D_\infty$ replaced by the  non-denegenerate limit $W_\infty$ of the martingale $W_n$.  By construction $W_\infty$ satisfies the almost sure invariance by random weighted mean equation  \begin{equation*}\label{fixedpoint} W_\infty= \sum_{|u|=1} \ee^{-V(u)}W_\infty(u), \end{equation*}  where $W_\infty(u)$ is the copy of $W_\infty$ built with the branching random walk associated the subtree of $\T$ rooted at $u$  (see \cite{KP,DL,Biggins76}), and it is worth recalling that the same holds for $D_\infty$ in the boundary case (see \cite{DL,K98,liu00}).

  % That $a \neq {3\over2}$ and can even take any value larger than or equal to 3 when the random walk $(S_n)$ has a finite variance [but has a non zero expectation $-\phi'(1-)$] reflects the fact that a different rule is involved in the behavior of typical  particles realizing $M_n$. 

  % this answers a question raised in \cite{AB}, where the multifractal analysis of the branching random walk in the first order phase transition is achieved. Since $W_\infty=0$ a.s. in case of a second order phase transition but then $Z_\infty$ is a non trivial non negative solution of \eqref{fixedpoint}, we find a remarkable parallel between the description of the asymptotic distribution of $M_n$ corresponding to the second and first order phase transition cases.
\medskip

We will state our assumptions in terms of the distribution of the i.i.d  increments $X_1,\ldots, X_n,\ldots$ of the random walk $(S_n)$ naturally associated with the branching random walk and assumed to be defined on a probability space  whose probability measure is   $\P$.  Denote by $\E$ the expectation with respect to $\P$ and set $X=X_1$. The law of $X$, denoted as $\P_X$,  is defined under \eqref{Hyp1}   by   
\begin{equation}
 \label{defS}
\int_{\mathbb R}f(x)\,\P_X(\mathrm{d}x): =\mathbb E\Big[\sum_{|u|=1}f(V(u)) \, \ee^{-V(u)} \Big]   ,
\end{equation}
for  any bounded measurable  function $f$.  Our first assumptions about $\P_X$ and expressed in terms of $X$ are the  following:   There exist  some constants $\gamma>3$, $  \alpha>1$,  a slowly  varying function $\ell$ and   some   $x_0 <0 $   such that  \begin{equation}
\label{Hyp2}
\mm:=\E[X] >0, \qquad  \E \Big[ (X^+)^\gamma\Big] < \infty, \qquad \P\left(  X  \le x \right) =\int_{-\infty}^x |y|^{-\alpha -1} \ell (y)  \mathrm{d} y , \qquad  \forall\,   x \le x_0,  \end{equation}

 \noindent with $y^+:= \max(y, 0)$ for any $y \in [-\infty, \infty)$.  The first property $m>0$ is just a restatement of $\phi(1-)<0$ whenever this derivative is defined.  The second and third properties imply in particular that $X$ is in the domain  of attraction of a stable law of index $\min(\alpha,2) $ (to fix ideas, let us mention that the boundary case considered in \cite{Aidekon13} correspond to  $\E[X] =0$ and $\E(X^2)<\infty$, as well as additional technical assumptions). One naturally gets a branching random walk leading to such an $X$ as follows: fix  a random variable $X$ obeying  \eqref{Hyp2} and assume in addition that $1<s=\E(\ee^X)=\int \ee^x \P_X(\mathrm{d}x)<\infty$ (in particular the second condition holds with all $\gamma>0$). Let  $(V_j)_{j\ge 1}$ be a sequence of random variables distributed according to $s^{-1}\ee^x \P_X(\mathrm{d}x)$,  $\nu$ a random integer independent of $(V_j)_{j\ge 1}$ and  such that $\E(\nu)=s$, and set $\L=\sum_{j=1}^\nu  \delta_{\{V_j\}}$. When $s$ is an integer, $\nu$ can be taken constant and equal to $s$, so that the branching random walk is built on the $s$-adic tree.

For  brevity, we extend the function $\ell$ to the whole $\r$,  by letting $\ell(x)= \ell(-x)$ for $x\ge |x_0|$  and $\ell(x) =1$ for any $x\in (x_0, |x_0|)$ [$|x_0|$ being large enough so that $\ell (x) >0$ for any $x \le x_0$].

\medskip

Under \eqref{Hyp1},  it is known that on the set ${ \bf  S}$ of  the survival of the system, $M_n \to \infty$ a.s. (see  Shi \cite{Shi12}). 

We have the following upper bound for the tightness of the minimum:

 \begin{Proposition}
\label{upbound} 
Under (\ref{Hyp1}) and  (\ref{Hyp2}), there exists some positive constant $K $ such that for all  $n\geq 2$ and  $x\geq 0$,
\begin{equation} \label{upbound12}
\mathbb P(M_n\leq \alpha_n -x)\leq  K  \ee^{-x},
\end{equation} where here and in the sequel, $$   \alpha_n:=  (\alpha+1)\log n -\log \ell(n).$$ 
\end{Proposition}

%An usual second moment computation and  Paley-Zygmund's inequality will give the  lower bound for the probability term in \eqref{upbound12}   if we assume (\ref{Hyp3}) and an  additional technical hypothesis (e.g. the forthcoming \eqref{tightness}). 

It is natural to study the convergence of  $M_n - \alpha_n$.  Before the presentation of the convergence in law under additional assumptions, let us say a few words on the norming constant $\alpha_n$.  For any $u \in \T\backslash\{\varnothing\}$, let ${\buildrel \leftarrow \over u}$ be the parent of $u$. Define \begin{equation} \label{def-Omega}
 \Delta V(u)  :=  V(u)- V( {\buildrel \leftarrow \over u}),  \qquad {\mathbb B}(u) :=   \Big\{ v:  v \neq u, {\buildrel \leftarrow \over v}= {\buildrel \leftarrow \over u}\Big\}. 
\end{equation}

\noindent For any $n\ge 1$ and  $ \vert u \vert = n$,  denote by $ \{ u_0:= \varnothing, u_1, ..., u_{n-1}, u_n =u\}$ the shortest path relating   the root $\varnothing$  to $u$ such that $\vert u_i\vert =i$ for any $0\le i \le n$.  

%their ancestors make a unique  huge  drop in the sense that 

It turns out  that the minimal position $M_n$ will be  reached  only by those particles $|u|=n$,  such that there is a unique $i\in [1, n]$  such that $\Delta V(u_i) < - n^{1+o(1)}$. Moreover, to make $V(u)=M_n$,  necessarily $i$ is near to $n$ and this (unique) large drop $\Delta V(u_i)$ will be of order $ -  n$,  which in view of the density function of $X$ in (\ref{Hyp2}), happens with probability of order $\ee^{-\alpha_n}$. This will yield the norming constant $\alpha_n$. However,  some particles $v \in {\mathbb B}(u_i)$ could also make a large drop in the sense that $\Delta V(v) < - n^{1+o(1)}$, moreover $v$ could also give some descendants which reach  $M_n$ in the $n$-th generation.   To get the convergence in law of $M_n - \alpha_n$,  we have to control this possibility of   simultaneous large drops in the same generation. This is why we need to introduce some extra conditions, stated below as   \eqref{decorr},  \eqref{tightness} and \eqref{tight+}.  We mention that these conditions hold for instance when   $\L = \sum_{i=1}^\nu \delta_{\{\xi_i\}}$ with $(\xi_i)$ i.i.d. and independent of $\nu$.

We also  need the following  integrability hypothesis, which combined with $\E(X)>0$,  is necessary and sufficient for  $W_\infty$ to not vanish almost surely \cite{Biggins76,KP,DL}: 
\begin{equation}
 \label{Hyp3}     \mathbb E \Big[ \big(\sum_{|u|=1} \ee^{-V(u)}\big)\,  \big(\log \sum_{|u|=1} \ee^{-V(u)} \big)^+\, \Big]<\infty,
\end{equation}

\noindent moreover $W_\infty>0$ on  ${ \bf  S}$.

\medskip

The main result of this paper is the following convergence in law:

\begin{theorem} \label{T:main}   
Assume (\ref{Hyp1}), (\ref{Hyp2}) and (\ref{Hyp3}), as well as    \eqref{decorr},  \eqref{tightness} and \eqref{tight+}.   Then for any $x\in \r$, 
\begin{equation} \label{mainconvergence}
\underset{n\to\infty}{\lim} \mathbb P(M_n\geq \alpha_n+x)=\mathbb E\left( \exp(- c_*\ee^x W_\infty)\right), 
\end{equation}
where $c_*>0$  is some constant   given in \eqref{defC*}.
 \end{theorem}

\begin{Remark}\label{R:unique} If almost surely $\#\{|u|=n: V(u)= M_n\}=1$ for any $n\ge1$, then we do not need the assumptions \eqref{decorr},  \eqref{tightness} and \eqref{tight+} in Theorem \ref{T:main}. 
\end{Remark}

In Theorem \ref{T:main}, the variety of possible behaviors obtained for $M_n$ comes for a part from the fact that a necessary and sufficient condition for the non degeneracy of $W_\infty$ is known, which makes it possible to choose $X$ with an infinite moment of order $\alpha$ with any $\alpha \in (1, \infty)$, while A{\"{\i}}d\'ekon \cite{Aidekon13}'s result assumes $\E(X^2)<\infty$ (and $\E(X)=0$), which with additional assumptions ensures that $D_\infty$ exists and is non degenerate; indeed, it is not known whether the assumption $\E(X^2)<\infty$ can be relaxed. 

%Interestingly, while this assumption  leads to the single $\frac{3}{2}\log n$ asymptotic behavior, in the context of Theorem~\ref{T:main} if one takes a constant function for  $\ell(x)$, then for each $\alpha\ge 2$ (as well as for each $\alpha\in (1,2)$) one finds  a different behavior, namely~$(\alpha+1) \log n$.

\medskip
Our result makes us conjecture that for $\beta>1$, the same  convergence result as in the boundary case holds for the Gibbs measures $\mu_{\beta,n}$ on $\{0,\ldots,s-1\}^{\N_+}$ if one replaces the critical Mandelbrot measure by the standard Mandelbrot measure, namely  the non degenerate measure which assigns mass $\ee^{-V(u)}W_\infty(u)$ to bond $u$. This would complete the parallel between the freezing phenomena observed under a second and a first order phase transition. The difference between these two situations can also be described at the critical temperature, and conditionally on non-extinction, as follows: under a second order phase transition, there exists a minimal supporting subtree $\T(0)$ for the free energy in the sense that the bounds of generation $n$ in $\T$ which mainly contribute to the free energy $F_n(1)$ are those $u$ of $\T(0)\cap\T_n$; moreover one observes the behavior, or singularity, $\frac{V(u)}{n}\approx 0=-\phi'(1)$ for the potential $V$ along $\partial\T(0)$, and  $\# \T(0)\cap\T_n\approx e^{o(n)}$. These properties are reminiscent  from the fact that in the infinite volume limit $\partial\T(0)$ is of Hausdorff dimension~0 and such that $\lim_{n\to\infty}\frac{1}{n}\log\sum_{|u|=n,\, [u]\cap \partial \T(0)\neq\emptyset}e^{-V(u)}=F(1)=0$, with $\lim_{n\to\infty}\frac{V(x_{|n})}{n}=0$ for all $x\in \partial \T(0)$, where $x_{|n}$ is the prefix of $x$ of length $n$ and $\partial\T$ is endowed with the standard ultrametric distance. Consequently, the free energy  concentrates on a single type of   singularity (see \cite{Moerters-Ortgiese,AB}). Under a first order phase transition, for all $\alpha\in [0,-\phi'(1)]$, there exists a subtree $\T(\alpha)$ of $\T$ such that $\#\T(\alpha)\cap\T_n \approx e^{n\alpha}$,  the bonds $u\in  \T(\alpha)\cap\T_n$ satisfy $\frac{V(u)}{n}\approx \alpha$, and they  substantially contribute to the free energy $F_n(1)$; in the infinite volume the  fractal sets $\partial\T(\alpha)$, $\alpha\in [0,-\phi'(1)]$, are of respective Hausdorff dimension $\alpha$, and such that $\lim_{n\to\infty}\frac{1}{n}\log\sum_{|u|=n,\, [u]\cap \partial T(\alpha)\neq\emptyset}e^{-V(u)}=F(1)=0$, and at each  $x\in \partial\T(\alpha)$ one observes the singularity $\lim_{n\to\infty}\frac{V(x_{|n})}{n}=\alpha$ (see \cite{AB} for more details). This can be interpreted as the coexistence of  uncountably many equilibrium states in the system at $\beta_c$.

\medskip

It is time to make  \eqref{decorr},  \eqref{tightness} and \eqref{tight+} explicit. To do so we need to introduce the probability measure $\q$ considered by Lyons~\cite{Lyons97} for general branching random walks (see also \cite{WW} for regular trees) and originally defined by Peyri\`ere~\cite{KP} for regular trees and in the case where $W_\infty$ is non degenerate ($\q$~is there defined as the skew product of the probability $\mathbb P$ and the Mandelbrot measure  $\mu$  \cite{KP} to study the Hausdorff dimension of $\mu$). 

\medskip

Denote by $({\cal F}_n, n\ge0)$ the natural filtration of the branching random walk.   The following  proposition is well-known:

% If the branching random walk starts from $V(\varnothing)=x$, then we denote its law by $\P_x$ (with $\P=\P_0$). 

%We still call $\mathbb{T}$ the genealogical tree of the process, so that $(w_n)_{n\in \N}$ is a ray of $\mathbb{T}$, which we will call the {\it spine}. This change of probability was also used in \cite{HuS09}. We refer to \cite{LPP95} for the case of the Galton-Watson tree, to \cite{CRo88} for the analogue for the branching Brownian motion, and to \cite{BKy04} for the spine decomposition in various types of branching.\nomenclature{$(w_n)_{n\in \N} $}{: spine of the branching random walk under $\q$}

\begin{Proposition}
\label{lyons}
Under (\ref{Hyp1}), on   the space  $\widehat \Omega$ of marked trees enlarged by an infinite distinguished  ray $(w_n, n \ge0)$, called spine, we may construct a probability measure $\q$  such that  

\medskip 

(i) for any $n\ge1$ and $|u|=n$, we have
\begin{equation}
\label{no(iii)}
\q \circ \pi^{-1} |_{\mathcal{F}_n}:=W_n\bullet\mathbb P|_{\mathcal{F}_n}, \qquad \q\left\{w_n=u\big| \pi^{-1}(\mathcal{F}_n)\right\}=\frac{\ee^{-V(u)}}{W_n},
\end{equation}
where  $\pi$ denotes the projection of $\widehat \Omega$ on $\Omega$;

\medskip 

(ii)  under $\q$,   $( \Delta V(w_n), \sum_{ v \in {\mathbb B}(w_n)} \delta_{\{ \Delta V( w_n)- \Delta V(v)\}})_{n \ge 1}$ is a sequence of i.i.d. random variables. Moreover,  the distribution of  $(V(w_n),n\geq 0)$  under $\q$  is the distribution of the random walk $(S_n,n\geq 0)$ under $\P$ defined above;

\medskip 

(iii) under $\q$, conditionally on ${\cal G}:= \sigma\{ u, \Delta V(u), {\buildrel \leftarrow \over u}=  w_j,  j \ge 0\}$,  the processes  $ \{ V(uv)- V(u), v \in \T\}$, for    $ u \in \cup_{j=1}^\infty {\mathbb B}( w_j), $ are i.i.d and   are distributed as $\{V(v), v \in \T\}$  under $\mathbb P$. 
\end{Proposition}

We refer the reader to \cite{ChauvinRouaultWakolbinger91, LyonsPemantlePeres95, Lyons97, BigginsKyprianou04, Shi12} for the detailed discussions on the change of measure and the proof of Proposition \ref{lyons}.  
  
\medskip

We denote by $\e_\q$ the expectation with respect to $\q$ and introduce the  first additional  hypothesis which we also  believe  necessary for the convergence of of $M_n - \alpha_n$: \begin{eqnarray} \label{decorr}
&&\nonumber \text{For any $f: \r \to \r_+$ measurable with compact support} \\
&&\lim_{z\to  -\infty}\mathbb E_\q \Big[ \ee^{- \sum_{ | v|= 1, v \neq w_1} f(   V(w_1)- V(v) ) } \, \big|\, V(w_1)= z\Big] \to \int \Xi(d \theta) \ee^{- \langle f,  \theta\rangle},
 \end{eqnarray}
 
 \noindent 
where $ \Xi$  is the distribution of some point process  on $\r\cup\{-\infty\}$ and we use the notation $\langle f, \theta\rangle:= \int_\r f(x) \theta(dx)$ for any $\theta \in {\cal M}$, the space of $\sigma$-finite measures on $\r\cup\{-\infty\}$.    For instance when   $\L = \sum_{i=1}^\nu \delta_{\{\xi_i\}}$ with $(\xi_i)$ i.i.d. and independent of $\nu$ it is easily seen that $ \Xi$ concentrates on $\delta_{\{-\infty\}} $. 

\medskip

The two other technical hypotheses are stated as follows:    
\begin{eqnarray}\label{tightness}  && \mbox{Under $\q$, as $z \to - \infty$, the   laws of $\#{\mathbb B}(w_1)$ conditionally on $\{V(w_1)= z\}$ are tight,}   \\
	\label{tight+}  &&  \lim_{\lambda \to \infty}  \limsup_{z \to -\infty}  \q \Big( \cup_{ v \in {\mathbb B}(w_k)} \{ \Delta V(w_k) - \Delta V(v) \ge \lambda\} \,   \big | \, \Delta V(w_k)= z \Big)    \to 0.
\end{eqnarray}

 It is  easy  to see that \eqref{tightness} and \eqref{tight+} are not very restrictive.  We shall explain the strategy of the proof of  Theorem~\ref{T:main} and Remark \ref{R:unique} in the next section.

 % \begin{Remark} \label{Remark1}    \begin{enumerate} \item  \eqref{Hyp1}  is a different from Durrett (1983) xxxxxx.   
%\item  Remark that we do not assume that the law of $X$ is diffuse. 
%\item xxxxxAdd a remark somewhere:   $\P(x < X \le x+ a)$ for any $a>0$ fixed and for $x \to -\infty$,   stronger than Bansaye and Vatutin etc xxxxxxx
% \item  Remark that \eqref{decorr} is stronger than the usual convergence in law of $\Xi_x $.
% \item 
 %When $\L = \sum_{i=1}^\nu \delta_{X_i}$ with $(X_i)$ i.i.d. and independent of $\nu$, \eqref{decorr} is satisfied with $ \Xi_\infty= \delta_{\{\infty\}}$ (??).
% \item Comparison with iid case.
%\end{enumerate}\end{Remark}

\section{Outline of the proof} \label{S2-0} 

 \begin{figure}[t]
 \centering
 \caption{}
 \label{Tux12}
 \includegraphics[interpolate=true,width=16.79cm,height=8.48cm]{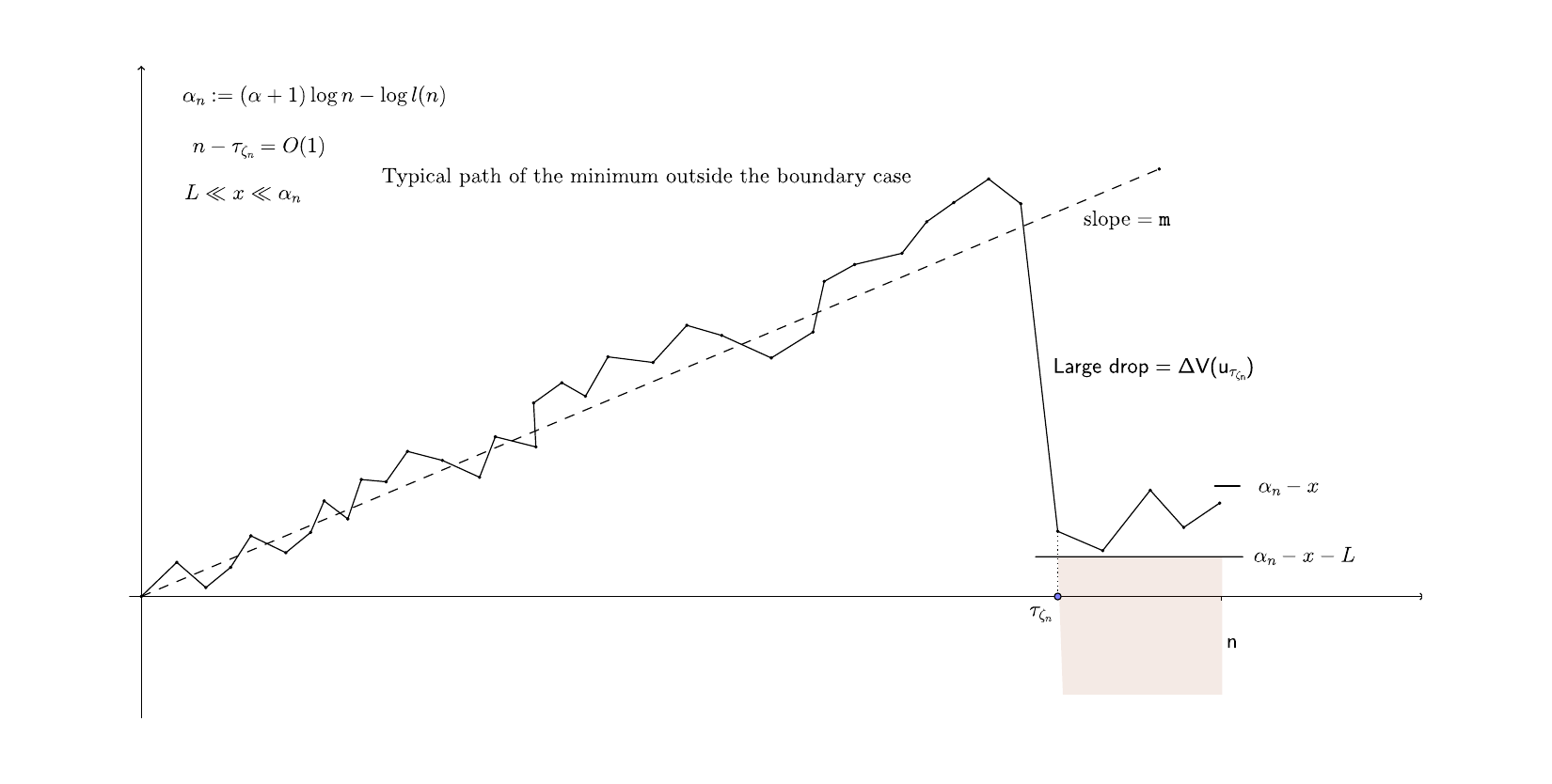} 
 \end{figure}

The main estimate  leading to Theorem \ref{T:main} is the following asymptotic tail for $M_n-\alpha_n$: 
 
\begin{Proposition}
\label{taildistrib}
Assume (\ref{Hyp1}), (\ref{Hyp2}),  (\ref{Hyp3}),     \eqref{decorr},  \eqref{tightness} and \eqref{tight+}.   For any $\varepsilon>0$, there exist  $A= A(\varepsilon)$ and an integer $n_0=n_0(\varepsilon)$ such that for all  $n>n_0$  and  $x\in [A,\frac{n}{\log n}]$, 
\begin{equation}
 \left| \mathbb P(M_n\leq \alpha_n  -  x)- c_* \ee^{-x}  \right|\leq \varepsilon  \, \ee^{-x}.
\end{equation}
\end{Proposition}
 
%\paragraph{Remark:} In comparison with the critical case $Cx\ee^{-x}$ is replaced by $C^* \ee^{-x}$.

\medskip

It turns out that the machinery developed by A{\"{\i}}d\'ekon in \cite{Aidekon13} is general enough to be adapted in the case considered in this paper. As a matter of fact, the proof  of  Proposition \ref{taildistrib} (of which Proposition~\ref{upbound} is one of the main ingredients) goes in the same spirit as  that of Proposition 1.3 in  A{\"{\i}}d{\'e}kon \cite{Aidekon13}, namely  the localization of the trajectory  of a particle  $u$    such that $V(u)= M_n$. The main difference is that, while in the boundary case such a trajectory typically corresponds to an excursion of length $n$, in our situation the trajectory $(V(u_j), 0\le j\le n)$  grows linearly until some generation $k$, near to $n$, where it makes a very large drop $\Delta V(u_k)$. 
% Such analysis will be presented in the proof of Proposition \ref{upbound} in Section \ref{S3}.  
To get Proposition \ref{taildistrib},  we shall prove that $n-k=O(1)$,  $\Delta V(u_k) =  -  (\mm +o(1))  n$ and  control the presence of several large drops in  the $k$-th generation  by using  the conditions \eqref{decorr},  \eqref{tightness} and \eqref{tight+}.

Specifically,  let us fix the threshold $\zeta_n:= {n\over (\log n)^3}$.   For any   $u\in \T$, 
let $\tau_{\zeta_n}^{(u)}$   be the first large drop in the path $\{V(u_i), 1\le i \le \vert u\vert\}$:
$$
\tau_{\zeta_n}^{(u)} : = \inf\{1\le   i \le |u|:   \Delta V(u_i)   <  -\zeta_n\},  $$

\noindent  with $\inf \emptyset:= \infty$.  Under the assumptions \eqref{Hyp1} and \eqref{Hyp2}, we analyze the particles leading to $M_n$ and obtain the following statement (see \eqref{bining1}):  Let $L$ and $T$ be large constants.   For all large $n  $ and for all $x>0$,  we have  $$\p\big( M_n\leq \alpha_n-x \big)  = \e\left[\frac{1}{\eta_n} \, \sum_{|u|=n} 1_{ \{M_n=V(u) \leq \alpha_n-x,\,  \min_{\tau_{\zeta_n}^{(u)}\le j \le n}V(u_j)  \geq \alpha_n-x-L ,\,    \tau_{\zeta_n}^{(u)} \in [n-T, n]\}}   \right] + o(1) \, \ee^{-x}  , $$

\noindent where $\eta_n:= \sum_{|u|=n} 1_{\{V(u)=M_n\}} $  and  $o(1) \to 0$  uniformly on $n$ and $x$,  as $L,  T \to \infty$.  

By the change of measure (cf. Proposition \ref{lyons}), the  above expectation is equal to $$\e_\q\left[\frac{1}{\eta_n} \, \ee^{V(w_n)}  1_{ \{M_n= V(w_n)\leq \alpha_n-x,\,  \min_{\tau_{\zeta_n}^{(w_n)} \le j \le n}V(w_j)  \geq \alpha_n-x-L ,\,     \tau_{\zeta_n}^{(w_n)} \in [n-T, n]\} }   \right] .$$

 The effects of simultaneous large drops are hidden in the number $\eta_n$, even if  at first sight this is not obvious.    Write $k:= \tau_{\zeta_n}^{(w_n)} \in [n-T, n]$. A crucial step in the localization of minimal particles, stated as Proposition \ref{clust}, says that under (\ref{Hyp1}), (\ref{Hyp2}) and (\ref{Hyp3}), for any $|u|=n$ such that $V(u)=M_n$, necessarily $u_{k-1}= w_{k-1}$, i.e.,  the trajectory of the particle $u$ and  the spine coincide at least until the generation $k-1$. Consequently, $\eta_n$ will only depend on the subtree rooted at $w_{k-1}$.    By the Markov property of the branching random walk under the  probability $\q$, we get that \begin{eqnarray*} &&  \p\big( M_n\leq \alpha_n-x \big) \\
	&&= \sum_{k=n-T}^n \ee^{\alpha_n -x} \,\mathbb E_\q\left[  1_{ \{ \tau_{\zeta_n} ^{(w_k)}=k,   V(w_k) \ge  y \}}    F^{(L)}_{n-k} \Big(  V(w_k)- y  , \sum_{v \in {\mathbb B}(w_k)} \delta_{\{ \Delta V(w_{k})- \Delta V(v)\}} \Big)  \right] + o(1) \, \ee^{-x}, \end{eqnarray*}

\noindent with some measurable function $F^{(L)}_j$ defined in \eqref{def-fjl} and $y:=\alpha_n-x-L$.  The next step will be an application of \eqref{decorr} [note that $\Delta V(w_{k}) \le - \zeta_n \to -\infty$] to get rid of the point measure $\sum_{v \in {\mathbb B}(w_k)} \delta_{\{ \Delta V(w_{k})- \Delta V(v)\}}$. Because of the compact support condition in \eqref{decorr}, we have to show that in this point measure, we can restrict ourselves to those $v \in {\mathbb B}(w_k)$ such that $ |\Delta V(w_{k})- \Delta V(v)|$ remains smaller than $\lambda$, with $\lambda>0$.  The hypotheses \eqref{tightness} and \eqref{tight+} are introduced to overcome this technical difficulty, as shown in the proof of Claim \ref{Claim-compact}.  Thus we get the truncated version of the above equality for $\p\big( M_n\leq \alpha_n-x \big)$  in \eqref{bmnx},  and an application of  \eqref{decorr} gives that  $$\mathbb P(M_n\leq \alpha_n-x)  =    \ee^{\alpha_n -x }\, \sum_{j=0}^{T-1} \E\left[  1_{ \{ \tau_{\zeta_n} =n-j, S_{n-j} \ge y\}} \,  G_j^{(\lambda,L)}(  S_{n-j} - y, \ X_{n-j})  \right] + o(1) \, \ee^{-x}, 
	$$ 

\noindent where  the measure function $G_j^{(\lambda, L)}(\cdot, \cdot)$ is defined in \eqref{def-gjl} and we have used the fact that under $\q$,  $(V(w_k), k\ge 0)$  is distributed  as the random walk $(S_k, k\ge0)$.   Finally, we apply a renewal result (Lemma \ref{renewal})  and get Proposition \ref{taildistrib}  by letting $\lambda, T, L  \to \infty$.  

Plainly,  if $\eta_n= 1$ a.s., then there is no effect coming from the possible simultaneous large drops and we get  Proposition \ref{taildistrib} without the assumptions  \eqref{decorr},  \eqref{tightness} and \eqref{tight+}, as stated in Remark \ref{R:unique}.

\medskip

Theorem \ref{T:main} follows from  Proposition \ref{taildistrib},  exactly as the main result in  A\"{i}dékon \cite{Aidekon13} follows from an analogous, though different, proposition (pp. 1405--1407). However, we give a proof for reader's convenience.

\medskip
{\noindent\it Proof of Theorem~\ref{T:main} as a consequence of Proposition \ref{taildistrib}:}  
For $B\ge 0$ define 
$$
\mathcal Z[B]=\{u\in\mathbb T: V(u)\ge B,\, V(u_k)<B,\, \forall k<|u|\}
$$
In the sense of \cite{K00} this is a very simple optional line and one has $\lim_{B\to\infty} \sum_{u\in \mathcal Z[B]} \ee^{-V(u)}=W_\infty$.

For $n\in\N_+$ and $0\le k\le n$, let $\Phi_{k,n}:x\ge 0\mapsto \mathbb P(M_{n-k}<\alpha_n-x)$.

Fix $x\in\mathbb R$ and $\varepsilon\in (0,c_*)$. Let $A(\varepsilon)$ and $n_0(\varepsilon)$ be defined as in Proposition~\ref{taildistrib}. Let $B>A(\varepsilon)+2|x|$ such that $(c_*+\varepsilon) \ee^{-B/2}<1$. Let $n_0\in\N_+$ such that $n_0\ge n_0(\varepsilon)$ and 
$$
\forall n\ge n_0, \ \p(\mathcal Y_{B,n})\ge 1-\varepsilon,
$$
where 
$$
 \mathcal Y_{B,n}=\{A(\varepsilon)\le V(u)-x\le \frac{n}{\log n},\, \forall u\in \mathcal Z[B]\}\cap \{\max\{|u|: u\in  \mathcal Z[B]\}\le n-n_0(\varepsilon)\}.
$$

Now for $n\ge n_0$ write 
$$
\mathbb P(M_n\ge \alpha_n+x)\ge \mathbb P(M_n\ge \alpha_n+x,\mathcal Y_{B,n})=\mathbb E\Big (\mathbf{1}_{\mathcal Y_{B,n}}\prod_{u\in\mathcal Z[B]}(1-\Phi_{|u|,n}(V(u)-x))\Big ),
$$
where we have used the conditional expectation along the stopping line. By construction we can apply  Proposition~\ref{taildistrib} to each term of the product and get
$$
\mathbb P(M_n\ge \alpha_n+x)\ge \mathbb E\Big (\mathbf{1}_{\mathcal Y_{B,n}}\prod_{u\in\mathcal Z[B]}(1-(c_*+\varepsilon)\ee^{x-V(u)})\Big )\ge\mathbb E\Big (\prod_{u\in\mathcal Z[B]}(1-(c_*+\varepsilon)\ee^{x-V(u)})\Big )-\p(\mathcal Y_{B,n}^c).
$$
 This yields 
 $$
 \liminf_{n\to\infty} \mathbb P(M_n\ge \alpha_n+x)\ge  \mathbb E\Big (\prod_{u\in\mathcal Z[B]}(1-(c_*+\varepsilon)\ee^{x-V(u)})\Big )-\varepsilon.
 $$
Moreover,  since $\max\{\ee^{-V(u)}:u\in\mathcal Z[B]\}$ tend a.s. to 0 as $B\to\infty$, we have  $\lim_{B\to\infty}\sum_{u\in\mathcal Z[B]}\log (1-(c_*+\varepsilon)\ee^{x-V(u)})=-(c_*+\varepsilon) \ee^x W_\infty$ hence by dominated convergence
$$
 \liminf_{n\to\infty} \mathbb P(M_n\ge \alpha_n+x)\ge\mathbb E(\exp (-(c_*+\varepsilon) \ee^x W_\infty))-\varepsilon,
 $$
and letting $\varepsilon$ tend to 0 yields the desired lower bound. To get the upper bound, write 
$$
\mathbb P(M_n\ge \alpha_n+x)\le \mathbb P(M_n\ge \alpha_n+x,\mathcal Y_{B,n})+\p(\mathcal Y_{B,n}^c).
$$
Following the same lines as above we get
$$
 \limsup_{n\to\infty} \mathbb P(M_n\ge \alpha_n+x)\le \mathbb E\Big (\prod_{u\in\mathcal Z[B]}(1-(c_*-\varepsilon)\ee^{x-V(u)})\Big )+\varepsilon,
 $$
and conclude as for the lower bound. 
\hfill$\Box$

\medskip

%In Section \ref{S2-0} we prove Theorem \ref{T:main} by admitting Proposition~\ref{taildistrib}.
The rest of the paper is organized as follows:  In Section \ref{S2}, we collect some preliminary estimates on  the one-dimensional random walk $(S_n)$, whereas we prove   Proposition \ref{upbound} in Section \ref{S3}.  Section  \ref{S:mainestimate} will use  \eqref{decorr},  \eqref{tightness} and \eqref{tight+} to  prove  Proposition \ref{taildistrib} by admitting a localization Lemma \ref{L:App1}. In  Section \ref{S:aux}, we give the proof of  Lemma~\ref{L:App1}.

Throughout the text,  we denote by $K$, $K'$ and $K''$  possibly with  several   subscripts,  some  positive constants whose values may change from one paragraph to another one.  We also wrote $f(n) \sim g(n)$ if $\lim_{n \to \infty} {f(n)\over g(n)} = 1$.

\section{Preliminaries on the one-dimensional random walk $(S_n)$.}\label{S2}

Recall that we considered in the introduction a sequence of i.i.d. real-valued random variables $(X_i)_{i\ge1}$ of common distribution that of $X$, and the random walk $(S_n)$ defined as  $S_n:= S_0+ X_1+...+ X_n$ for any $n \ge1$ with $S_0\in \r$.   Let  $ \overline S_n:= \max_{0\le k \le n} S_k$ and  $\underline S_n:= \min_{0\le k \le n} S_k$.   For  $x \in \r$,    denote  the  distribution  of  $(S_n)$ by $\P_x$ if $S_0=x$ and    $\P= \P_0$. We state some known facts as lemmas:

\begin{Lemma}[\cite{DenisovDiekerShneer08}, pp. 1950, Lemma 2.1] 
\label{L:1}  Let $(S_n)$ be a one-dimensional random walk  satisfying  that $\E[  \vert S_1\vert^b ] < \infty $   for $b>1$.     Let $\mm:= \E(S_1)$.  There exists a constant $K= K_b>0$ such that for all $n\geq 1,\, y\geq  n^{\max({1\over b}, {1\over2})}$ and $x>0$,
\begin{eqnarray}
\P(S_n- \mm n  \le -x,\, \underset{1\leq i\leq n}{\min} X_i\geq -y)  &\leq & K  \ee^{-\frac{x}{y}},  \label{L:1a}\\ 
	\P(|S_n- \mm n| \ge x,\, \underset{1\leq i\leq n}{\max} |X_i|\leq y) &\leq & K  \ee^{-\frac{x}{y}}. \label{L:1b}
\end{eqnarray}
\end{Lemma}

\begin{Lemma}[Gut \cite{Gut09},  Theorem 6.2, pp. 93] \label{L:Gut}
Let $S$ be a one-dimensional random walk with positive mean $\mm$ starting from $0$.  Let 
$$ R(x):= \sum_{n=0}^\infty \P \Big( \overline S_n \le x \Big), \qquad x\ge0.$$
Then $$ \lim_{x \to \infty} { R(x)\over x}= \mm.$$
\end{Lemma}

%As a consequence, for any $n\ge 1$ and $x\ge 0$, $  n \P( \overline S_n \le x) \le \sum_{i=1}^n \P( \overline S_i \le x) \le R(x)$, hence $$  \P\Big( \overline S_n \le x\Big)  \le {R(x)\over n}.$$ 

%\begin{Lemma}[Nagaev \cite{Nagaev82}] Assume that $\e[ \vert S_1\vert ]  < \infty$ and $\p\Big( S_1 < x \Big) = {1\over \alpha} \vert x\vert^{-\alpha} \ell(\vert x\vert )$, as $x\to - \infty$, with some $\alpha>1$. Then  for any fixed $\gamma>0$, uniformly on $x\ge \gamma n$, $$ \P\Big( S_n - \E[S_n] \le - x \Big) \sim n\, \P\Big( S_1 \le -x \Big). $$ Moreover there exists some $n_0\ge1$ such that for all $n\ge n_0$ and $x \ge \gamma n$,  $$ \P\Big( S_n - \E[S_n] \le - x \Big)  \le 2 n \, \P\Big( S_1 \le -x \Big). $$   \end{Lemma}

%\begin{Lemma}[Borovkov and Borovkov \cite{Borovkov08}, Theorems 8.2.3,  8.2.18 and 8.3.4,  Bertoin and Doney \cite{BertoinDoney96}]  \label{L:BBo08}    Assume \eqref{Hyp2}. We have that  for any fixed $x \ge0$, $$ \P\Big( \overline S_n \le x \Big) \sim e^D U(x) \P(S_1 < -  \mm  n),$$ with $D:= ?, U(x):=?$ Moreover for any $x \to \infty$ such that $x \le {\mm\over2} n$, we have $$ \P \Big( \overline S_n \le x \Big) \sim {x\over \mm}  \, \P\Big( S_1 < x- \mm n\Big).$$ By Stone \cite{Stone65}'s local limit theorem (and the boundedness of the stable density, see Vatutin and  Wachtel \cite{VatutinWachtel09} for $\ell_1$),  \end{Lemma}

\begin{Lemma}[Stone \cite{Stone65}]  \label{L:BBo08}    Assume \eqref{Hyp2}.  There exists  some slowly varying function  $\ell_1$ such that for all $x \in \r$, $\hslash  >0$ and $n \ge 1$,   $$ \P\Big( S_n \in [x, x + \hslash ]\Big) \le   \hslash \,n^{- \max(1/\alpha, 1/2)}   \ell_1(n) .   $$      \end{Lemma}

 We mention that up to  a multiplicative  constant,  $\ell_1$  only depends on   the truncated second moment of $S_1$, see Vatutin and  Wachtel \cite{VatutinWachtel09}. In particular, 
 if $\alpha>2$, we may choose $\ell_1\equiv  K$ for some  positive  constant large enough.

 Let us  introduce the drops in the random walk  $(S_n)$:   for   $\zeta>0$,  define 
\begin{eqnarray} \label{drop1}  \tau_{\zeta} &:= &\inf\{j\geq 1: \, X_j  <  -\zeta\},  \\
	 \label{drop2} \tau_\zeta^{(2)} & :=& \inf\{j> \tau_{\zeta}: \, X_j  <  -\zeta\} , \end{eqnarray} the first    and the second    drop  of size $\zeta$.    We shall consider  $\zeta \in [{\zeta_n \over4}, 4 \zeta_n]$ with    \begin{equation}\label{zetan} \zeta_n:=\frac{n}{(\log  n)^3}, \qquad n\ge 2. \end{equation}

\begin{Lemma}\label{3urd} Assume \eqref{Hyp2}.   There exists  some constant $K>0$   such that for all  $n \geq 2$, $-\infty < y\leq {\mm\over2} n $,  $$ \P\left( S_n-y \in [0,1]\right)  \le  K\, n^{-\alpha} \ell(n).$$
\end{Lemma}

{\noindent \it Proof of Lemma \ref{3urd}:}  It is enough to consider large $n$.  Observe that 
\begin{eqnarray}
\nonumber \P\left( S_n-y \in [0,1]\right)  &\le &  \P(S_n-y\leq 1,\,  \tau_{\zeta_n} > n )+  \P(S_n -y\in [0,1],\tau_{\zeta_n}^{(2)}  \le n)  \nonumber \\
	\qquad \qquad && + \sum_{i=1}^n\P(S_n-y\in[0, 1],\, \tau_{\zeta_n}=i,  \tau_{\zeta_n}^{(2)} > n )   \nonumber
\\  &=:&  A_{\eqref{3urdn}} + B_{\eqref{3urdn}}+C_{\eqref{3urdn}}.  \label{3urdn}
 \end{eqnarray}

 For any  $y\leq {\mm\over2} n$,    \begin{equation}\label{A3urdn} A_{\eqref{3urdn}} \le \P( S_n - \mm n \le  1- {\mm\over2} n,  \tau_{\zeta_n} > n ) \le K\ee^{-   ({\mm\over2} n-1)/\zeta_n} \le \ee^{- {\mm\over3} (\log n)^3}, \end{equation}

 \noindent  where we have applied Lemma \ref{L:1} to get the second inequality in \eqref{A3urdn}.  For $B_{\eqref{3urdn}}$, we deduce from Lemma \ref{L:BBo08}   that  
 \begin{eqnarray*} B_{\eqref{3urdn}} &= & \sum_{i=1}^{n-1} \sum_{j=i+1}^n \P \Big( \tau_{\zeta_n}=i, \tau_{\zeta_n}^{(2)}=j, S_n - y \in [0, 1]\Big) \\
 	 &\le& \sum_{i=1}^{n-1} \sum_{j=i+1}^n \P \Big( \tau_{\zeta_n}=i, \tau_{\zeta_n}^{(2)}=j\Big)   \, (n-j+1)^{- \max(1/\alpha,1/2)} \ell_1(n-j+1) \\
	&\le&  n \times \zeta_n^{-2\alpha} \ell(\zeta_n)^2\,  n^{ 1 - \max(1/\alpha,1/2)}  \max_{1\le k \le n} \ell_1(k) = o(    n^{-  \alpha  } ), 
 \end{eqnarray*}

 \noindent since $2-\alpha-\max(1/\alpha,1/2)<0$.  Finally for all $n \ge 2$,    let  \begin{equation} \label{def-ein}
{E_i^{(n)}}:=\Big\{ | S_n-X_i - \mm (n-1) | \leq \frac{n}{\log n}  \Big\}  , \qquad 1\le i \le n.
\end{equation}

 Observe that for    any ${\zeta_n \over4} \le \zeta \le  4\zeta_n$\footnote{We consider $\zeta$ instead of $\zeta_n$ for the use of \eqref{1probab1} in the proof of Lemma \ref{P:1}; Moreover, by exchangeability, the probability $\P(    \tau_{\zeta}=i  ,\, \tau_{\zeta}^{(2)} >n ,\, ({E_i^{(n)}})^c)$ does not depend on $i$.} and $1\le i \le n$,  \begin{eqnarray}
\nonumber && \P\Big(    \tau_{\zeta}=i  ,\, \tau_{\zeta}^{(2)} >n ,\, ({E_i^{(n)}})^c\Big)
\\ &&=  \P\Big(    \tau_{\zeta}=n    ,\, ({E_n ^{(n)}})^c\Big) \nonumber
   \\
\nonumber &&=  \P\left( X_n  <   -\zeta\right)  \P\Big(\min_{1\le j\leq n-1} X_j\geq -\zeta,\, |S_{n-1}-\mm (n-1)| > \frac{ n}{\log n} \Big)
 \\
 &&\leq \zeta^{-\alpha} \, \ell(\zeta) \left ( \P\big(\max_{1\le j\leq n-1}X_j\geq \zeta \big) + \P\Big(\max_{1\le j\leq n-1}|X_j|\leq \zeta,\, |S_{n-1}-\mm (n-1)| > \frac{n}{\log n} \Big)\right )  \nonumber\\
	&&\le    	\zeta^{-\alpha} \, \ell(\zeta) 	\Big( n  \zeta^{-\gamma} \E[(X^+)^\gamma] +  K \, \ee^{-\frac{n }{\zeta \log n}} \Big) \qquad \mbox{(by using 	(\ref{Hyp2}) and (\ref{L:1b}))} 
	\nonumber \\ && \le     n^{-(\alpha+\gamma-1)} \ell_2(n),
	\label{1probab1}
\end{eqnarray}

\noindent  with some  slowly varying function $\ell_2$.   Using  the exchangeability, \begin{eqnarray*}
C_{\eqref{3urdn}} &=& n\,  \P \Big( X_n < - \zeta_n,  S_n - y \in [0, 1], \tau_{\zeta_n}=n\Big)  \\
	& \le&  n^{-(\alpha+\gamma-2)} \ell_2(n) + n \, \E\Big[ 1_{  E^{(n)}_n} \P\big( X_n + s - y \in [0, 1], X_n < -\zeta\big)\big\vert_{s=S_{n-1}}\Big],
	\end{eqnarray*}

\noindent by using the independence of $X_n$ and $S_{n-1}$. Notice that  \begin{equation} \label{Karamata} 
  \sup_{x \le - {\mm\over3} n }\, \frac{\ell(   x   )}{|x|^{\alpha+1}}  \le (1+o(1))   \ell(  n) \big({\mm\over3} n \big)^{-\alpha-1} \le K\,  \ee^{-\alpha_n}, 
\end{equation}
 
 \noindent by  using Karamata's representation for the slowly varying function $\ell$.  On $E^{(n)}_n$,    $y-S_{n-1} \le -{\mm\over2} n +m +  {n\over \log n}\le - {\mm\over3} n-1$. It follows from \eqref{Hyp2}  and \eqref{Karamata} that  on $E^{(n)}_n$, uniformly for $s= S_{n-1}$,  $\P\big( X_n + s - y \in [0, 1]\big)  \le    (1+o(1))   \ell(  n) \big({\mm\over3} n \big)^{-\alpha-1} $, which  implies that  for all large $n$, $C_{\eqref{3urdn}} \le  \big({\mm\over3}   \big)^{-\alpha-1}(1+o(1))   n^{-\alpha} \ell(  n).$    Lemma \ref{3urd} follows from $\eqref{3urdn}$. \qed

\medskip

Recall that $\alpha_n=  (\alpha+1)\log n -\log \ell(n)$.

\begin{Lemma}
\label{P:1}   Assume \eqref{Hyp2}.    There exist  $K >0$ and some slowly varying function $\ell_3\ge1$  such that   for all  large $n\geq n_0$, $\forall\, \zeta \in [{\zeta_n \over4}, 4 \zeta_n]$,  $ a\le {n \over \log n}$,    we have that     \begin{equation} \label{pasprov2bis-sum}   \P\left(S_n-y\in [a,a+1 ],\,  \min_{  \tau_\zeta \le j \le n} S_j   \geq y,\,   \tau_{\zeta}^{(2)}\leq n \right)  \leq       R(a+1  ) \, n ^{1-\max(1/\alpha,1/2)  -2\alpha} \,  \ell_3(n) , \qquad  \forall y \in \r,
\end{equation} 
whereas  for all  $- \infty <   y <  {\mm\over 2} n $,  \begin{equation}
\label{pasprov-sum}   \P\left( S_n-y\in [a,a+1 ],\,    \min_{  \tau_\zeta \le j \le n} S_j  \geq y,\,   \tau_{\zeta} \le n <  \tau_{\zeta}^{(2)}   \right)  \leq   K \,      R( a+1 ) \,     \ee^{-\alpha_n}     .\end{equation}
Moreover, for any $\varepsilon>0$, there exists some $\lambda=\lambda(\varepsilon  )>0$ such that for all  $- \infty <   y <  {\mm\over 2} n $,  \begin{equation}
\label{tight2}   \P\left( S_n-y\in [a,a+1 ],\,  \min_{  \tau_\zeta \le j \le n} S_j  \geq y,\,   | S_n - S_{\tau_\zeta}| >\lambda,  \tau_\zeta \le n    \right)  \leq    \varepsilon\,     \ee^{-\alpha_n}     .\end{equation}
\end{Lemma}

\medskip
The similar results hold if we replace the interval $[a, a+1]$ by $[a, a+\hslash]$ with an arbitrary positive  constant $\hslash$. 
\medskip

\noindent{\it Proof of Lemma \ref{P:1}.}    We shall prove that for   any $1\le i <  n$, 
\begin{equation} \label{pasprov2bis}   \P\Big(S_n-y\in [a,a+1 ],\, \min_{  \tau_\zeta \le j \le n} S_j   \geq y,\, \tau_{\zeta}=i,\,  \tau_{\zeta}^{(2)}\leq n \Big)  \leq      R(a+1 \, ) \, i^{-\max(1/\alpha,1/2)} \,   n^{-2\alpha} \,  \ell_3(n) , \qquad  \forall y \in \r, 
\end{equation} 

\noindent whereas  for all $- \infty <   y <  {\mm\over 2} n $,  \begin{equation}
\label{pasprov}   \P\Big( S_n-y\in [a,a+1 ],\,  \min_{  \tau_\zeta \le j \le n} S_j  \geq y,\, \tau_{\zeta}=i,\,    \tau_{\zeta}^{(2)} > n \Big)  \leq   K \,      \P\Big(\overline{S}_{n-i+1}\leq a+1 \,  \Big)   \, \ee^{-\alpha_n}    +  n^{-(\alpha+\gamma-1)} \ell_3(n).\end{equation}

Clearly, up to a   multiplicative constant,  \eqref{pasprov2bis-sum}  and  \eqref{pasprov-sum} follow   from \eqref{pasprov2bis} and  \eqref{pasprov} by taking the  sum over $i \in 1, 2, ...,  n-1$.

Let us denote by   $ \P_{(\ref{pasprov2bis})}(i)$ the probability term  in (\ref{pasprov2bis}). By considering the time-reversal random walk $(S_n-S_{n-k},\, 0\leq k\leq n) {\buildrel  (d)  \over = }(S_k,\, 0\leq k\leq n)$,  we get  that  for any $1\le i \le n-1$,  $(S_n,   \min_{  \tau_\zeta \le j \le n} S_j  ,   \{\tau_{\zeta}=i< \tau^{(2)}_\zeta \le n\}  )$ has the same distribution as $(S_n,   S_n - \overline S_{\sigma_n},  \{\sigma_n= n-i+1> \tau_{\zeta} \})$, where $\sigma_n:= \max\{k\in [1,n],\, X_k <  -\zeta\}$ (with the usual convention $\max \emptyset:=0$).  It follows that 
\begin{eqnarray*}
    \P_{(\ref{pasprov2bis})}(i) &=&  \P\left( S_n-y\in [a,a+1 ],\, \overline{S}_{n-i+1}\leq S_n- y,\,  \sigma_n=n-i+1  > \tau_{\zeta}  \right) 
	\\ &\leq &  \P\left( S_n-y\in [a,a+1 ],\, \overline{S}_{n-i}\leq a+1 \, ,\,  X_{n-i+1} < -\zeta,\,  \tau_{\zeta} < n-i +1\right)  
	\\&=& \E\left[1_{ \{   X_{n-i+1} < -\zeta,\,  \tau_{\zeta} < n-i+1 ,  \overline{S}_{n-i+1}\leq a+1 \, \}  } \P_{S_{n-i+1}}\left(  S_{i-1}-y\in [a,a+1 ]\right) \right],  
	\end{eqnarray*}

\noindent by the Markov property at $n-i+1$.  Set $g(i)=\sup_{z\in\r} \P_z\big(  S_{i}-y\in [a,a+1] \big)$.  We have 
	\begin{eqnarray}
	 \P_{(\ref{pasprov2bis})}(i)  &\le &      g(i-1) \,  \P\Big(   X_{n-i+1} < -\zeta,\,  \tau_{\zeta} < n-i ,  \overline{S}_{n-i}\leq a+1 \, \Big)   \nonumber \\
	\nonumber &\leq &     g(i-1) \,   \sum_{1\leq j<n-i } \P\left( X_{n-i} < -\zeta,  X_j <  -\zeta,\, \overline{S}_{j-1}\leq a+1 \,    \right)  	\\
\nonumber &=  &    g(i-1) \,  \P\left( X  < -\zeta\right)^2  \,     \sum_{1\leq j <  n- i}  \P\left( \overline{S}_{ {j-1 }}\leq a+1 \, \right)
\\
\label{YUE2} &\leq &      g(i-1) \,  \zeta^{-2\alpha} \, \ell(\zeta)^2\, R(a+1 \, ),
\end{eqnarray}

\noindent for all large $n$. According to Stone's local limit theorem (Lemma \ref{L:BBo08}), there exists a constant $C>0$ such that for $i\ge 2$ one has $g(i-1)\le   i^{-\max(1/\alpha,1/2)} \ell_1(i) $, and since $g(0)\le 1$, $C$ can be chosen so that $g(0)\le C \ell_1 (1)$.  This   yields \eqref{pasprov2bis} as we shall choose   $$\ell_3(n): = \max( \ell_2(n),    4^{2\alpha}\,    (\log n)^{6\alpha} \,  \max_{1\le i \le n, {\zeta_n\over4} \le \zeta \le 4 \zeta_n} C\,\ell_1(i)  \ell(\zeta)^2),$$ where $\ell_2(n)$ is the slowly varying function appeared in \eqref{1probab1}.

%
%With  Lemma \ref{L:BBo08}, by combining (\ref{YUE1}) and (\ref{YUE2}) we obtain
%\begin{eqnarray}
%\nonumber A_{(\ref{proba1})}+ C_{(\ref{proba1})} & \leq &C  \ee^{-\frac{n-x}{n}(\log n)^3}+   n^{1-\frac{1}{\alpha}} \zeta_n^{-2\alpha} l(\zeta_n)^2 R(a+1 \, )
%\\
%\label{YUE3} &\leq &(a+1 \, )(1+1 \, ) o(\mm^{-(\alpha+1)} \, \ee^{-\alpha_n} ).
%\end{eqnarray}
%
%
%It remains to estimate $B_{(\ref{proba1})}$ which corresponds to the case where we have only one big jump below $-\zeta_n$. As
%\begin{eqnarray}
%\nonumber B_{(\ref{proba1})}&=& \sum_{i=1}^n   \P(S_n-y\in [a,a+1 \, ],\,   \underline{S}_{[i,n]}\geq y,\,  X_i\leq -\zeta_n ,\, \min_{j\neq i}X_j \geq -\zeta_n),
%\end{eqnarray}
%we apply (\ref{pasprov}) then Lemma \ref{L:Gut} and deduce that
%\begin{eqnarray}
%\nonumber   B_{(\ref{proba1})}&\leq & \mm^{-(\alpha+1)} \, \ee^{-\alpha_n} \sum_{i=1}^n \P\left( \overline{S}_{n-i}\leq a+1 \, \right)\leq  C(a+1 \, )1 \,  \mm^{-(\alpha+1)} \, \ee^{-\alpha_n} .
%\end{eqnarray}
%\\

  To prove    (\ref{pasprov}), we first establish an inequality implying that when $S_n=o(n)$, with a big probability there is a unique large drop  $X_{\tau_{\zeta}}  $  before $n$  which is of order of magnitude $-\mm n$.   Recall \eqref{def-ein} for  the definition of ${E_i^{(n)}}$.     Define for any $i \in [1, n]$, 
\begin{equation}
\label{cugon}
 \P_{(\ref{cugon})}(i):=  \P\Big(S_n-y\in [a,a+1 ],\,   \min_{  i  \le j \le n} S_j  \geq y,\,  \tau_{\zeta}=i  ,\, \tau_{\zeta}^{(2)} >n ,\, {E_i^{(n)}}\Big) .
\end{equation}

 In view of  \eqref{1probab1},  (\ref{pasprov}) will follow if we can prove   that  \begin{equation} \label{cugon2} \P_{(\ref{cugon})}(i) \le   K\,    \P\left(\overline{S}_{n-i+1}\leq a+1 \,  \right)  \, \ee^{-\alpha_n} . \end{equation}

 \noindent By conditioning on $\sigma\{X_j, 1\le j \le n, j \neq i\}$,  we have that  \begin{eqnarray}  \P_{(\ref{cugon})}(i)  &\le& \P \Big(  S_n -  \min_{  i  \le j \le n} S_j  \le  a +1 \, ,  S_n-y\in [a,a+1 ],  {E_i^{(n)}}\Big)  \nonumber \\
 	&=&  \E\Big[ 1_{ \{ S_n -  \min_{ i  \le j \le n} S_j  \le  a +1 \, , {E_i^{(n)}}\}} \P \big( X_i +t - y \in [a, a+1]\big)\big\vert_{t= S_n- X_i}\Big]. \label{adad1}
 	\end{eqnarray}

  On ${E_i^{(n)}}$, $ \vert t-  \mm (n-1)\vert \le  {n \over \log n}$, then $z\equiv a+1 \,  + y - t \le - {\mm\over3} n$ for all large $n\ge  n_0 $ and uniformly for all $a \le {n\over \log n}$ and $ y < {\mm\over2} n$, hence it follows from  \eqref{Hyp2} and \eqref{Karamata} that   $$   \P_{(\ref{cugon})}(i) \le   K\,  \ee^{-\alpha_n}\,   \P \Big(S_n -  \min_{ i  \le j \le n} S_j  \le  a +1 \, , {E_i^{(n)}}   \Big)    ,$$

\noindent which yields \eqref{cugon2} by using the fact that $ \P(S_n -  \min_{ i  \le j \le n} S_j  \le  a +1 ) = \P ( \overline S_{n-i+1} \le a + 1 ).$ 
This completes   the proof of \eqref{pasprov}.

Remark that in \eqref{adad1}, if we replace the event $\{S_n -  \min_{  i  \le j \le n} S_j  \le  a +1\}$ by $\{ | S_n - S_i|  >\lambda\}$ with   $\lambda>0$, then   for any $i \in [1, n]$, \begin{equation} \label{adad2} 
\P \Big(  | S_n -  S_i| > \lambda \, ,  S_n-y\in [a,a+1 ] ,  {E_i^{(n)}}\Big) \le  K\,  \ee^{-\alpha_n} \, \P(  | S_n - S_i|  >\lambda) . 
\end{equation}

   Denote by $\P_{\eqref{tight2}}$ the probability term in \eqref{tight2}.   Notice that by     \eqref{pasprov2bis-sum},  the probability that the event in \eqref{tight2} holds together with $\{\tau_\zeta^{(2)} \le n\}$ is bounded by  $ R(a+1) n ^{1-\max(1/\alpha,1/2)  -2\alpha} \,  \ell_3(n) \le  {\varepsilon\over 4}\,    \ee^{-\alpha_n}$ for all large $n\ge n_0(\varepsilon)$.  On the other hand,  we deduce from  \eqref{pasprov} that for some  large but fixed  integer $k=k(\varepsilon,  a)$,  \begin{eqnarray*}  &&  \P\left( S_n-y\in [a,a+1 ],\, \min_{\tau_{\zeta}\le j\le n}S_j\geq y,\,  \tau_{\zeta} <  n-k ,\,    \tau_{\zeta}^{(2)} > n \right)  
  	\\ &\leq&   K \,      \sum_{j=k}^\infty  \P\left(\overline{S}_j \leq a+1 \,  \right)   \, \ee^{-\alpha_n}    +  n^{1 -(\alpha+\gamma-1)} \ell_3(n)
   \\ &\le & {\varepsilon\over 4}\,    \ee^{-\alpha_n} ,
   \end{eqnarray*}
   
   \noindent  for all $n \ge n_1(\varepsilon)$ [recalling that $\gamma>3$].   Therefore   \begin{eqnarray*}
\P_{\eqref{tight2}} & \le  &   {\varepsilon\over 2}\,    \ee^{-\alpha_n}  +  \P\left( S_n-y\in [a,a+1 ],   | S_n - S_{\tau_\zeta}| >\lambda,  n-k \le   \tau_\zeta\le n  <  \tau_{\zeta}^{(2)}   \right) 
	\\ &\le&  {\varepsilon\over 2}\,    \ee^{-\alpha_n}  +  (k+1)\,  n^{-(\alpha+\gamma-1)} \ell_2(n)+   \sum_{i=n-k}^n  \P\left( S_n-y\in [a,a+1 ],   | S_n - S_{\tau_\zeta}| >\lambda,      \tau_\zeta= i  , E_{i}^{(n)}  \right) , 
	\end{eqnarray*}

   \noindent by applying \eqref{1probab1} to $i= n, n-1, ..., n-k$.  Since $\gamma>3$,   $(k+1)\,  n^{-(\alpha+\gamma-1)} \ell_2(n) \le {\varepsilon\over 4}\,    \ee^{-\alpha_n}$,  which in view of   \eqref{adad2} imply that  $\P_{\eqref{tight2}}  \le {3\varepsilon\over 4 }\,    \ee^{-\alpha_n} +   K\, \sum_{j=0}^k \P\big( |S_j| > \lambda\big) \, \ee^{-\alpha_n} \le \varepsilon \ee^{-\alpha_n}$, if we choose some $\lambda=\lambda(k, \varepsilon)$ large enough. This proves \eqref{tight2} and completes the proof of Lemma \ref{P:1}.   \qed

\medskip
We present  a renewal result associated to the random walk $(S_n)_{n\geq 0}$.
\begin{Lemma}
\label{renewal}
Under  (\ref{Hyp2}). Let $G:\r_+ \times \r  \to \r$ be a  measurable  function such  that for some $b >1$  and some positive constant $K>0$, 
\begin{equation}\label{control}
 \sup_{z \in \r} | G(x, z) | \le    K \,  (1+x)^{-b}, \qquad \forall\, x \ge0.
 \end{equation} 
 Assume furthermore that for any $x \in \r_+$,      $\lim_{z \to -\infty}  G(x, z) $ exists, and denote it by  $G_*(x)$. Then   
\begin{equation}
\label{eqrenewal}
\lim_{n\to\infty} \ee^{\alpha_n}  \E\left[ 1_{\{ \tau_\zeta=n,  \, S_n \geq y \}} G(S_n-y,  X_n )\, \right]=  \mm^{-(\alpha+1)}   \int_{0}^\infty  G_*(x)\, \mathrm{d} x,  
\end{equation}  uniformly on $ | y|  \leq \frac{n}{\log n}$ and ${\zeta_n\over4} \le \zeta \le 4 \zeta_n$. 
\end{Lemma}

\noindent{\it Proof of Lemma (\ref{renewal}).}  Without loss of generality we may assume that $G$ takes nonnegative values.  Let $\varepsilon>0$ be small.  Let $E^{(n)}_n:= \{ |S_{n-1}-\mm  (n-1) |\leq \frac{n}{\log n}\}$ as in \eqref{def-ein}.  By \eqref{1probab1}, $$ \P\big( \tau_\zeta= n, (E^{(n)}_n)^c\big) \le n^{-(\alpha+\gamma-1)} \ell_2(n) \le \varepsilon   \, \ee^{-\alpha_n} , $$

\noindent for all large $n$.  Let us denote by $\E_{\eqref{eqrenewal}}$ the expectation term in \eqref{eqrenewal}.  Then  \begin{equation} \ee^{\alpha_n}  \E_{\eqref{eqrenewal}} =  \ee^{\alpha_n}  \E\left[  1_{ \{ \tau_\zeta =n , S_n \geq y \} \cap E^{(n)}_n} \,   G(S_n-y,  X_n )  \right]   + O(\varepsilon)  
 . \label{eqre1}\end{equation}

To deal with the above expectation term, we  distinguish two situations according to the value of $S_n -y$:      Clearly,  \begin{eqnarray}  \ee^{\alpha_n}  \E\left[  1_{ \{ \tau_\zeta =n , S_n - y\ge  {n\over \log n}   \} \cap E^{(n)}_n} \,   G(S_n-y,  X_n )  \right]   &  \le&   \ee^{\alpha_n} K\,  ( 1+ {n\over \log n}     ) ^{- b} \P ( X_n < - \zeta )   \nonumber
	\\ &= &  K\,  ( 1+{n\over \log n} )^{- b}  \ee^{\alpha_n}  \int_{- \infty}^{-\zeta} |x|^{-\alpha-1} \ell (x)\, \mathrm{d} x \nonumber
	\\ &\le&  \varepsilon ,  \label{eqre2} \end{eqnarray}

\noindent     uniformly on $\zeta\in [{\zeta_n \over 4}, 4 \zeta_n]$ since $b>1$.   If $0\le  S_n - y  < {n\over \log n}  $, then on the event $E^{(n)}_n$,  $X_n= S_n - S_{n-1} $ satisfies that $| X_n  +  \mm  (n-1)| \le 3  {n \over \log n}  $   uniformly on $|y|\le {n \over \log n}$, hence $\ee^{\alpha_n}  |X_n|^{-\alpha-1} \ell (X_n)=\mm ^{-(\alpha+1)} +o(1)$.  Moreover, since for $n$ large enough on $E_n^{(n)}$ we have $X_n\le -\zeta$ so that $\tau_\zeta\le n$, and it is easily seen that conditionally on $\tau_\zeta\le n$, the probability that $\tau^{(2)}_\zeta>n$ tends to $1$ as $n$ tends to $\infty$, we can write $1_{ \{ \tau_\zeta =n , 0 \le S_n - y< {n\over \log n}   \} \cap E^{(n)}_n}=1_{E^{(n)}_n}(1-\gamma_n)$, with $\lim_{n\to\infty}\gamma_n=0$ uniformly on $|y|\le {n \over \log n}$.

Therefore for all large $n$, 
 \begin{eqnarray}  &&  \ee^{\alpha_n} \E\left[  1_{ \{ \tau_\zeta =n , 0 \le S_n - y< {n\over \log n}   \} \cap E^{(n)}_n} \,   G(S_n-y,  X_n )  \right]   \nonumber \\
	&= &   \E\Big[ 1_{E^{(n)}_n} \,  \int_{\substack{ \{  0 \le S_{n-1}+z  - y< {n\over \log n}   \} \\\{ |z + \mm  (n-1) | \le 3{n\over\log n}\}}  } \,   G(S_{n-1}+z-y,  z) ( \ee^{\alpha_n}  |z|^{-\alpha-1} \ell (z))\, \mathrm{d} z \Big] +R_n,
 \end{eqnarray}	
with 
$$
R_n= - \ee^{\alpha_n} \E\left[ \gamma_n \ee^{\alpha_n}\cdot  1_{ \{ \tau_\zeta =n , 0 \le S_n - y< {n\over \log n}   \} \cap E^{(n)}_n} \,   G(S_n-y,  X_n )  \right]. 
$$
This yields
 \begin{eqnarray}  &&  \ee^{\alpha_n} \E\left[  1_{ \{ \tau_\zeta =n , 0 \le S_n - y< {n\over \log n}   \} \cap E^{(n)}_n} \,   G(S_n-y,  X_n )  \right]   \nonumber \\
	 \nonumber
	\\ & =&   (\mm ^{-(\alpha+1)} +o(1))\,  \E\Big[ 1_{ \{ \tau_\zeta =n\}\cap E^{(n)}_n} \,  \int 1_{ \{  0 \le S_{n-1}+z  - y< {n\over \log n}   \}   } \,   G(S_{n-1}+z-y,  z)    \, \mathrm{d} z \Big]  +R_n \nonumber \\
	&=& (\mm ^{-(\alpha+1)} +o(1))\,  \E\Big[  1_{ E^{(n)}_n} \int_0^{{n\over\log(n)}}G(x,  x+y - S_{n-1})    \, \mathrm{d} x \Big]  +R_n \nonumber \\
	&=:&  (\mm ^{-(\alpha+1)}+o(1))   \, \E_{\eqref{eqre3}} +R_n. \label{eqre3}
 \end{eqnarray}

Notice that for any fixed $ x \in \r_+$ and  $|y| \le {n\over \log n}$, $ 1_{ E^{(n)}_n}\,       G(x,  x+y - S_{n-1})$ converges a.s.   to $  G_*(x)  $    as $n \to \infty$. Indeed  $S_{n-1}$ tends linearly to  $-\infty$ and $ 1_{E^{(n)}_n}$ converges to 1 a.s. by the Kolmogorov-Marcinkiewicz-Zygmund law of large numbers. 
 
 It then follows from \eqref{control}, the dominated convergence theorem, and the fact that $y - S_{n-1}$ tends a.s. uniformly to 0 on  $|y| \le {n\over \log n}$ that $\E_{\eqref{eqre3}} \to  \int_0^\infty      G_*(x)    \mathrm{d} x$,  uniformly on $|y| \le {n\over \log n}$.  Then, bounding the function integrated in $R_n$ by $\ee^{\alpha_n}1_{E^{(n)}_n} \,   G(S_n-y,  X_n )$ and using bounded convergence theorem we get $R_n\to 0$, still uniformly on $|y| \le {n\over \log n}$. In view of \eqref{eqre1},  \eqref{eqre2} and \eqref{eqre3}, this yields the desired conclusion.  \hfill$\Box$

\section{Proof of Proposition \ref{upbound}} \label{S3}

At first let us fix some notations which will be used throughout the rest of this paper: For $ \vert u \vert = n$, we write  $[\varnothing, u]\equiv \{ u_0:= \varnothing, u_1, ..., u_{n-1}, u_n =u\}$ the shortest path  from the root $\varnothing$ to $u$ such that $\vert u_i\vert =i$ for any $0\le i \le n$.  For any $u, v \in \T$, we use the partial order  $ u < v$ if $u$ is an ancestor of $v$ and  $u \le v$ if $u < v$ or $u=v$. By the standard  words-representation in a tree, $u < v$  if and only if the word $v$ is a concatenation of the word $u$ with some word $s$, namely $v = us$ with
$|s| \ge 1$.  Denote by $\T^{(u)}:= \{ v : u \le v \}$ the subtree rooted at $u$ and by $\T_n:= \{v: |v|=n\}$ the set of vertices at generation $n$ for any integer $n$.  Let   ${\buildrel \leftarrow \over v}$  be  the parent of $v$ for any $v \neq \varnothing$.

%For $a\in \r$, we denote by $\P_a$ the probability distribution associated to the branching random walk starting from $a$, and $\E_a$ the corresponding expectation.  
\medskip

The following  so-called many-to-one formula (\ref{manytoone}) can be obtained   as a consequence of the spinal decomposition (see Proposition \ref{lyons}): 
Under (\ref{Hyp1}), for any $n\geq1$ and any measurable function $g:\r^n\to[0,+\infty)$,
\begin{equation}
\label{manytoone}
\mathbb E\Big[\underset{|u|=n}{\sum}g(V(u_1),...,V(u_n))\Big]=\E\Big[\ee^{S_n}g(S_1,...,S_n)\Big].
\end{equation}

The proof of  Proposition \ref{upbound} will be based on the forthcoming three lemmas. The first one is a well-known fact in the studies of branching random walk: 

\begin{Lemma}\label{L:expmin} Assume \eqref{Hyp1}. We have that  \begin{equation}\label{trebas}  \mathbb P \Big( \exists u\in \T,\, V(u)\leq - x \Big) \le \ee^{-x },  \qquad \forall \, x >0. \end{equation}
\end{Lemma}

  Lemma \ref{L:expmin}  follows  from  a simple application of  \eqref{manytoone} if one considers the first generation  $k$ such that for    some  $u \in \T_k$,   $V(u) \le -x$, see e.g. Shi \cite{Shi12} for details.

%{\noindent \it Proof of Lemma \ref{L:expmin}.} We give the proof only  for the completeness:  By considering the first generation $k\ge1$ for which there exists some $\vert z \vert=k$ such that $V(z) \le -x$, we get that  \begin{eqnarray} \nonumber \P\left(\exists z\in \mathbb{T},\, V(z)\leq - x \right)&\leq& \sum_{k=1}^{\infty}\E\Big[ \sum_{|z|=k}  1_{\{ V(z)\leq -x ,\, V(z_i)>-x,\, \forall i\leq k \}}\Big] \\ \nonumber &\leq & \sum_{k=1}^{\infty}  \E\Big[\ee^{S_k}  1_{\{S_k\leq -  x ,\, S_i> -x,\, \forall i\leq k \}}\Big] \\ \nonumber &\leq & \ee^{-x} \sum_{k=1}^{\infty}  \P\left(   S_k\leq -x,\, S_i > -x,\, \forall i\leq k \right)  \\ \label{trebas}&\leq &\ee^{-x}\P\left(\exists k\in \N,\,  S_k\leq -x\right)=o(\ee^{-x}), \end{eqnarray}  where in the second inequality we used the many-ton-one formula \eqref{manytoone}.  \qed

%To prove Proposition \ref{upbound} we will continuously use Lemma \ref{P:1} with $1 \, $ fixed equal to $1$. 

To state the second lemma, we need to introduce some notations similar to that in \eqref{drop1} and \eqref{drop2}: Recall that $\zeta_n:= {n\over (\log n)^3}$.   For any   $u\in \T$, 
let $\tau_{\zeta_n}^{(u)}$ and $\tau_{\zeta_n}^{(2, u)}$ be the first and the second large drop in the path $\{V(u_i), 1\le i \le \vert u\vert\}$:
\begin{eqnarray}
\tau_{\zeta_n}^{(u)} &: = &\inf\{  i \in [1, \vert u \vert  ]:   V(u_i)-V(u_{i-1})  <  -\zeta_n\},  \label{drop-u1} \\
	 \tau_{\zeta_n}^{(2, u)}&:=& \inf\{  i \in (\tau_{\zeta_n}^{(u)}, \vert u \vert ]:   V(u_i)-V(u_{i-1}) <  -\zeta_n\}, \label{drop-u2}
\end{eqnarray}  with $\inf \emptyset:= \infty$.   Recall    that $\alpha_n= (\alpha+1) \log n - \log \ell(n)$.   Our second lemma says that for those $u$ such that $V(u) \le \alpha_n -x$,  necessarily there is an  unique large drop before $\vert u\vert$:

\begin{Lemma}
\label{forcedjump}
Assume (\ref{Hyp1})  and (\ref{Hyp2}).  For any $\varepsilon >0$ there exists $n_0(\varepsilon)>0$ such that for any $n\geq n_0(\varepsilon)$ and all $ x\geq 0$,
\begin{eqnarray}
 \p\big( \exists  u\in \mathbb{T}_n,\, V(u)\leq \alpha_n-x,\,   \tau_{\zeta_n}^{(u)}>n\big) &\le& \varepsilon \,    \ee^{-x}  ,   \label{eqforcejump1} \\
 \p\big( \exists u\in \mathbb{T}_n,\,  V(u)\leq \alpha_n-x,\, \min_{\tau_{\zeta_n}^{(u)}  \le    j \le n}V(u_j)  \ge -x- \alpha_n ,\,  \tau_{\zeta_n}^{(2, u)}\leq  n\big)  &\le& \varepsilon\,    \ee^{-x}  . \label{eqforcejump2}
\end{eqnarray}
Consequently for any $x >0$, \begin{equation} 
 \p\big( \exists u\in \mathbb{T}_n,\,  V(u)\leq \alpha_n-x,\,   \tau_{\zeta_n}^{(2, u)}\leq  n\big)  \le  \varepsilon\,    \ee^{-x}   \label{eqforcejump2-new}
 \end{equation}
\end{Lemma}

 \medskip
 We may replace in \eqref{eqforcejump2} $\min_{\tau_{\zeta_n}^{(u)}  \le    j \le n}V(u_j)  \ge -x- \alpha_n  $ by $\min_{\tau_{\zeta_n}^{(u)}  \le    j \le n}V(u_j)  \ge -x-  n^{b}$ with any constant $b \in (0,  \alpha+  {1\over\alpha} -2)$.

\medskip

\noindent{\it Proof of Lemma \ref{forcedjump}.}  By the many-to-one formula  \eqref{manytoone} and using the notations \eqref{drop1} and \eqref{drop2},  the probability term in  \eqref{eqforcejump1} is less than  \begin{eqnarray*}  \mathbb E\Big(\sum_{|u|=n}  1_{\{ V(u)\leq \alpha_n-x,\tau_{\zeta_n}^{(u)}>n\}}\Big) &=  &\E\left(\ee^{S_n}  1_{\{ S_n\leq \alpha_n-x,\tau_{\zeta_n} >n\}}\right)
	\\ &\leq & \ee^{- x+ \alpha_n}\P\left(S_n\leq \alpha_n-x,\tau_{\zeta_n} >n \right) 
	\\& \le &    \ee^{- x+ \alpha_n} \P\Big(S_n - \mm n \leq \alpha_n - \mm n , \min_{1\le i \le n} X_i \ge - \zeta_n \Big)  
	\\ &\le&      \ee^{-{\mm\over 2}  (\log n)^3}  \ee^{-x} ,  \end{eqnarray*}    for all large $n \ge n_1$ and 
where we have used Lemma \ref{L:1}  for the last inequality.   This proves \eqref{eqforcejump1}.

\medskip

Let us  denote by  $\p_{\eqref{eqforcejump2} }$   the      probability term in \eqref{eqforcejump2}. Then  
\begin{eqnarray*}
     \p_{\eqref{eqforcejump2} } &\le &  \e\left[\sum_{|u|=n}  1_{\{   V(u)\leq \alpha_n-x,\,    \min_{\tau_{\zeta_n}^{(u)} \le j \le n} V(u_j) \ge -x- \alpha_n,  \tau_{\zeta_n}^{(2, u)}\leq  n\}}\right]  
   	\\ &  =  &  \E\left[\ee^{S_n}  1_{\{ S_n\leq \alpha_n-x,\, \min_{\tau_{\zeta_n}\le j \le n}S_j \geq -x- \alpha_n,\,   \tau_{\zeta_n}^{(2)}\leq  n\}}\right]
	\\ & \leq &  \sum_{k= 1}^{  \lceil 2 \alpha_n \rceil +1}  \ee^{k - x - \alpha_n }\P\left(   S_n+   \alpha_n+x  \in [k-1, k),\,    \min_{\tau_{\zeta_n}\le j \le n}S_j  \geq -x-\alpha_n ,\, \tau_{\zeta_n}^{(2)}\leq  n \right).
\end{eqnarray*}

By applying (\ref{pasprov2bis-sum})  with    $y\equiv  -x- \alpha_n$,   we get that  [$R$ is a nondecreasing function] for any $1\le k \le   \lceil 2 \alpha_n \rceil +1$,  $$\P\left(   S_n+   \alpha_n+x  \in [k-1, k),\,    \min_{\tau_{\zeta_n}\le j \le n}S_j  \geq -x-\alpha_n ,\, \tau_{\zeta_n}^{(2)}\leq  n \right)  \le  n^{ 1- 1/\alpha - 2\alpha} \ell_3(n)\, R( 2\alpha_n +2 ),$$

\noindent which implies that $$ \p_{\eqref{eqforcejump2} } \le  \ee^{-x -\alpha_n} \, \ee^{2 \alpha_n+1} \,   n^{ 1- 1/\alpha - 2\alpha} \ell_3(n)\,  R( 2\alpha_n+2 )= \ee^{-x}\; n^{ 2- \alpha- 1/\alpha} \,  \ell_4(n)\,$$ with some slowly varying  function $\ell_4$. Since for $\alpha>1$, $2- \alpha- 1/\alpha < 0$,  \eqref{eqforcejump2} follows.   

\medskip

Finally,   we deduce from \eqref{eqforcejump2} and Lemma \ref{L:expmin} that  the probability term in  \eqref{eqforcejump2-new} is less than $ \varepsilon \ee^{-x} + \P\big( \exists u \in \T: V(u) < - x- \alpha_n\big) \le \varepsilon \ee^{-x} + \ee^{-\alpha_n -x} \le 2 \varepsilon \ee^{-x}$ yielding   \eqref{eqforcejump2}. \hfill$\Box$

\medskip
Below is the third and the last  lemma that we need in the proof of Proposition \ref{upbound}:

\begin{Lemma}
\label{mauvaistraj}
Assume (\ref{Hyp1})  and (\ref{Hyp2}). There exist $K,\, c_4>0$ such that  for any $n$ and $L_0\in \N^*$ large enough, and for any  $x\ge 0$ and $ L\in [L_0, (2+\alpha)\log n]$, 
\begin{equation}
\label{intermL}
\p\big( \exists u\in \mathbb{T}_n,\, V(u)\leq \alpha_n-x,\, \min_{\tau_{\zeta_n}^{(u)}\le j \le n}V(u_j) -( \alpha_n-x)\in [-L,-L+1],\, \tau_{\zeta_n}^{(u)}\leq n<\tau_{\zeta_n}^{(2, u)}\big)\leq  K\, \ee^{-c_4L}\ee^{-x}.
\end{equation}
Consequently there exists some constant $c_2>0$ such that for any $L\ge   L_0$, \begin{equation}
\label{intermL2}
\p\big( \exists u\in \mathbb{T}_n,\, V(u)\leq \alpha_n-x,\, \min_{\tau_{\zeta_n}^{(u)}\le j \le n}V(u_j)  \le   \alpha_n-x -L,\, \tau_{\zeta_n}^{(u)}\leq n<\tau_{\zeta_n}^{(2, u)}\big)\leq  K\, \ee^{-c_2 L}\ee^{-x} .
\end{equation}
\end{Lemma}

%\begin{Lemma}
%\label{bontraj}
%Under (\ref{Hyp1}), (\ref{Hyp2}) and (\ref{Hyp2}), for any $n,\, L_0\in \N^*$ large enough and for any  $x \geq 0,\, L\in [L_0, (2+\alpha)\log n]$,
%\begin{equation}
%\label{eqbontraj}
%\P\big( \exists z\in \mathbb{T}_n,\, V(z)\leq \alpha_n-x,\,  \underline{V}_{[ \tau_{\zeta_n}^{(z)},n]}(z)\geq \alpha_n-x-L,\,  \tau_{\zeta_n}^{(z)}\leq n<\tau_{\zeta_n}^{(2, z)}\big)\leq L^2\ee^{-x}.
%\end{equation}
%\end{Lemma}

\noindent{\it Proof of Lemma \ref{mauvaistraj}.} Let $ \mathbb P_{(\ref{intermL})}$ the probability term  in (\ref{intermL}).  Pick up a constant  $\beta \in (0,\frac{1}{4(2+\alpha)})$. Notice that $L<(2+\alpha)\log n$ implies  $      \ee^{\beta L}     \le n^{ \frac{1}{4}}$.  For notational simplification, we write in this proof $$ y \equiv y(n, x, L):= \alpha_n - x - L$$
(notice that $y<\mm n/2$ if $n$ is large enough). 

For any $u\in \mathbb{T}_n$ satisfying the condition in the probability term in (\ref{intermL}), there exists $p \in [\tau_{\zeta_n}^{(u)}, n] $ such that $V(u_p)\in   [y, y+1]$. Then $\tau_{\zeta_n}^{(u)}=\tau_{\zeta_n}^{(u_p)}$, and

\begin{eqnarray}
\nonumber \p_{(\ref{intermL})}&\leq  &\sum_{p=1}^n   \p\Big( \exists u\in \mathbb{T}_n,\,  \min_{\tau_{\zeta_n}^{(u)}\le j \le n}V(u_j)   \geq y,\, V(u_p) -y \in [0,1], V(u)\leq y +L,\, \tau_{\zeta_n}^{(u)} \le p,\,   \tau_{\zeta_n}^{(2, u)} >n  \Big)
\\
\label{poipo}& \leq & \sum_{p=1}^{n-      \lfloor  \ee^{\beta L} \rfloor     }    A_{(\ref{poipo})}(p)  +  \sum_{p=n-       \lfloor  \ee^{\beta L}     \rfloor}^n B_{(\ref{poipo})}(p),
\end{eqnarray}

\noindent
with
\begin{eqnarray*}
 A_{(\ref{poipo})}(p)&:= &\e\left[\sum_{|u|=n} 1_{\{ \min_{ \tau_{\zeta_n}^{(u)}\le j \le n}V(u_j)   \ge y ,\,  {V}(u_p) -y \in [0, 1],\, V(u)\leq  y+L,\, \tau_{\zeta_n}^{(u)}\le p ,\,  \tau_{\zeta_n}^{(2, u)} >n  \}}  \right],
 \\
 B_{(\ref{poipo})}(p)&:=&  \e\left[\sum_{|v|=p}  1_{\{    \min_{ \tau_{\zeta_n}^{(v)}\le j \le p}V(v_j)  \ge y ,\, {V}(v) -  y\in [0, 1],\, \tau_{\zeta_n}^{(v)} \le  p<\tau_{\zeta_n}^{(2, v)}  \}}   \right],
\end{eqnarray*}

\noindent where   the sum of the expectation term of $ B_{(\ref{poipo})}(p)$ is obtained by considering $v=u_p$ satisfying $V(u_p) \in [y, y+1]$.  We omitted the dependence on $n$ in both $A_{(\ref{poipo})}(p)$ and $B_{(\ref{poipo})}(p)$.  By using (\ref{manytoone}), we have
\begin{eqnarray}
   B_{(\ref{poipo})}(p)&=&  \E\Big[\ee^{ S_p}   1_{\{ \min_{\tau_{\zeta_n}\le j\le  p}S_j  \ge y ,\, S_p - y \in [0, 1], \, \tau_{\zeta_n} \le  p<\tau_{\zeta_n}^{(2)} \}}  \Big]
	 \nonumber\\ &\le& \ee^{y+1} \, \P\Big(   \min_{\tau_{\zeta_n}\le j\le  p}S_j   \ge y ,\, S_p - y \in [0, 1], \, \tau_{\zeta_n}  \le  p<\tau_{\zeta_n}^{(2)}\Big) 
	 \nonumber\\ &\leq &    K'\, \ee^{y} \,     \ee^{-\alpha_p} ,
\label{Bpoipo}  
\end{eqnarray}

\noindent where the last inequality follows from     (\ref{pasprov-sum}) by remarking that ${\zeta_p\over4} \le \zeta_n \le 4\zeta_p$ for any     $ p \in [ n-  \lfloor  \ee^{\beta L} \rfloor, n]$; moreover   $\ee^{-\alpha_p} \sim   \ee^{-\alpha_n} $,  so  for all large $n $,    $$  \sum_{p= n-\lfloor  \ee^{\beta L} \rfloor }^n B_{(\ref{poipo})}(p)  \le  2 \,   K'    \ee^{y + \beta L}     \, \ee^{-\alpha_n}   \le   K\, \ee^{-x - L/2}. $$

It remains to estimate  $A_{(\ref{poipo})}(p)$.   By applying    (\ref{manytoone}), \begin{eqnarray}
\nonumber  A_{(\ref{poipo})}(p)&=&  \E\Big[\ee^{S_n}  1_{\{    \min_{\tau_{\zeta_n}\le j\le  n}S_j   \ge   y ,\,  S_p - y \in [0, 1]  ,\,  S_n \leq y+L ,\,  \tau_{\zeta_n} \le p ,\, \tau_{\zeta_n}^{(2)} > n \}}  \Big]
\\
\nonumber &\leq &  \ee^{y+L} \P \Big(     \min_{\tau_{\zeta_n}\le j\le  n}S_j  \geq  y,\, S_p- y \in [0,1],\, S_n-y  \in [0,L],\,   \tau_{\zeta_n}\le p   ,\, \tau_{\zeta_n}^{(2)} >n  \Big) . 
\end{eqnarray}

\noindent By applying    the Markov property at time $p$, we see that the above probability term is  equal to $$ \E\Big[   1_{\{   \min_{\tau_{\zeta_n}\le j\le  p}S_j  \geq  y,\, S_p- y \in [0,1],   \tau_{\zeta_n} \le p < \tau_{\zeta_n}^{(2)} \}} \P_{S_p} \big( \underline{S}_{n-p } \geq  y , S_{n-p} - y \in [0, L], \tau_{\zeta_n} > n-p\big)   \Big]. 
 $$  
 
 \noindent For any $z\equiv S_p \in [y, y+1]$,  $ \P_{z} \big( \underline{S}_{n-p } \geq  y , S_{n-p} - y \in [0, L], \tau_{\zeta_n} > n-p\big)  \le  \P  \big( \underline{S}_{n-p } \geq  -1  , S_{n-p}   \in [-1, L+1], \tau_{\zeta_n} > n-p\big) $. Recalling  that $y =\alpha_n -x - L$,   we get  
\begin{equation}\label{poipo12} A_{(\ref{poipo})}(p) \le  \ee^{\alpha_n -x } \,  I_{\eqref{poipo12}} \, J_{\eqref{poipo12}} , \end{equation}
 
 \noindent with \begin{eqnarray*}
I_{\eqref{poipo12}} &:=& \P\big(    \min_{\tau_{\zeta_n}\le j\le  p}S_j  \geq  y,\, S_p- y \in [0,1],    \tau_{\zeta_n} \le p <  \tau_{\zeta_n}^{(2)}     \big)  , \\
J_{\eqref{poipo12}} &:=& \P  \big( \underline{S}_{n-p } \geq  -1  , S_{n-p}   \in [-1,L+ 1], \tau_{\zeta_n} > n-p\big). 
\end{eqnarray*}

 For $1\le p <  \lfloor {n\over 4}\rfloor $, we  apply  Lemma \ref{L:1}  to see that  $ J_{\eqref{poipo12}}   \le   \P  \big(  S_{n-p}  \le L+ 1 , \tau_{\zeta_n} > n-p\big) \le K \ee^{ - ( \frac{3\mm n}{4} -(L+1))/\zeta_n}  $, hence for     $1\le p <  \lfloor {n\over 4}\rfloor $, 
 \begin{equation}
\label{Bpoipo2} A_{(\ref{poipo})}(p)\le  K\,      \ee^{\alpha_n -x } \,     \ee^{ - ( \frac{3\mm n}{4}-(L+1))/\zeta_n}  \leq     \ee^{-{\mm \over2} (\log n)^3}   \ee^{-x} .
\end{equation}

For $\lfloor \frac{n}{4} \rfloor \le   p \le n-\lfloor  \ee^{\beta L} \rfloor$,  we   apply   (\ref{pasprov-sum}) for $I_{\eqref{poipo12}}$  (with $y\equiv \alpha_n-L-x \le  \alpha_n  \le {\mm \over2} p$ and $\zeta=\zeta_n \in [{\zeta_p\over4}, 4\zeta_p]$),   and $L+1$ times Lemma \ref{3urd}  for $ J_{\eqref{poipo12}} $ (recall that $n-p\geq   \lfloor  \ee^{\beta L}   \rfloor\geq \ee^{\beta L_0}$ large and thus $L\le {\mm\over 2} (n-p)$) we get
$$  A_{(\ref{poipo})}(p)\leq  K'\,  \ee^{\alpha_n -x }  \,  \ee^{-\alpha_n}  \,    (L+1)\frac{\ell(n-p)}{(n-p)^{\alpha}}   \leq  K{''}\, \ee^{-x  } (L+1) \frac{\ell(n-p)}{(n-p)^{\alpha}} , $$

\noindent which together with    (\ref{Bpoipo2})  yield  that 
\begin{eqnarray*} \sum_{p=1}^{n-\lfloor  \ee^{\beta L} \rfloor} A_{(\ref{poipo})}(p) & \le &    \ee^{-x} \,  n\,  \ee^{-{\mm \over2} (\log n)^3} +  K^{''} \ee^{-x} (L+1)  \sum_{p=\lfloor \frac{n}{4} \rfloor }^{n-\lfloor  \ee^{\beta L} \rfloor}   \frac{\ell(n-p)}{(n-p)^{\alpha}} 
	\\ &\le& K\,  \ee^{-x} \ee^{- \frac{(\alpha-1)}{2} \beta L}  ,
\end{eqnarray*}

\noindent proving (\ref{intermL}).

It remains to prove \eqref{intermL2}.  Let $c_3:= {1\over 2(\alpha+2)}$. Remark that if $L \ge {\alpha_n\over 1-c_3}$, then $\alpha_n -L \le - c_3 L$ and it follows   from Lemma \ref{L:expmin} that  $$ 
\p\big( \exists u\in \mathbb{T}_n,  \min_{\tau_{\zeta_n}^{(u)}\le j \le n}V(u_j)  \le   \alpha_n-x -L,\, \tau_{\zeta_n}^{(u)}\leq n \big) \le \ee^{- c_3 L -x}.$$

Therefore it is enough  to treat the case $L_0\le  L < {\alpha_n\over 1-c_3}$. As $n \ge n_0$, $ {\alpha_n\over 1-c_3} \le (2+\alpha) \log n$. Then the probability term in  \eqref{intermL2} is less than 

 \begin{eqnarray*}&&\p\big ( \exists u\in \mathbb{T}_n,\, V(u)\leq \alpha_n-x,\, \min_{\tau_{\zeta_n}^{(u)}\le j \le n}V(u_j) -( \alpha_n-x)\le -c_3L,\, \tau_{\zeta_n}^{(u)}\leq n<\tau_{\zeta_n}^{(2, u)}\big) \\&&+ \sum_{k=L}^{\alpha_n+c_3 L} \p\big( \exists u\in \mathbb{T}_n,\, V(u)\leq \alpha_n-x,\, \min_{\tau_{\zeta_n}^{(u)}\le j \le n}V(u_j) -( \alpha_n-x)\in [-k,-k+1),\, \tau_{\zeta_n}^{(u)}\leq n<\tau_{\zeta_n}^{(2, u)}\big) \\
&& \le\ee^{-c_3 L -x} + \sum_{k=L}^{\alpha_n/(1-c_3)}\ee^{-c_4 k-x} \\ &&\le 
\ee^{-c_3 L -x} + K' \ee^{-c_4 L -x}.
\end{eqnarray*}

\noindent  We get  \eqref{intermL2} by choosing $c_2: =\min(c_3, c_4)$.  \hfill$\Box$

 \medskip
Now we can tackle the {\it Proof of Proposition \ref{upbound}.}  We fix an arbitrary integer  $L\in[ L_0, (2+\alpha) \log n]$ (as in Lemma \ref{mauvaistraj}) and consider large  $n$. Then 
\begin{equation}
  \p\left(M_n\leq \alpha_n-x\right)\leq  \p_{\eqref{mnx1234}}^{(1)} + \p_{\eqref{mnx1234}}^{(2)}+\p_{\eqref{mnx1234}}^{(3)}+\p_{\eqref{mnx1234}}^{(4)} , \label{mnx1234}
\end{equation}
with
\begin{eqnarray*}
\p_{\eqref{mnx1234}}^{(1)} &:=& \p\left(\exists u\in \mathbb{T},\, V(u)\leq -x\right)
\\
\p_{\eqref{mnx1234}}^{(2)}&:=& \p\big( \exists u\in \mathbb{T}_n,\,  V(u)\leq \alpha_n-x,\, \min_{\tau_{\zeta_n}^{(u)}\le j \le n}V(u_j) >  -x,\,  \tau_{\zeta_n}^{(2, u)}\leq  n\big) \\
	&& \qquad\qquad +\p\big( \exists u\in \mathbb{T}_n,\, V(u)\leq \alpha_n-x,\,   \tau_{\zeta_n}^{(u)}>n\big)
\\
\p_{\eqref{mnx1234}}^{(3)}&:=& \sum_{k=L}^{  \alpha_n   +1    } \p\big( \exists u\in \mathbb{T}_n,\, V(u)\leq \alpha_n-x,\, \min_{\tau_{\zeta_n}^{(u)}\le j \le n}V(u_j)  -( \alpha_n-x)\in [-k,-k+1),\,\tau_{\zeta_n}^{(u)}\leq n<\tau_{\zeta_n}^{(2, u)} \big)
\\
\p_{\eqref{mnx1234}}^{(4)}&:=& \p\big( \exists u\in \mathbb{T}_n,\, V(u)\leq \alpha_n-x,\, \min_{\tau_{\zeta_n}^{(u)}\le j \le n}V(u_j)  \geq \alpha_n-x-L,\, \tau_{\zeta_n}^{(u)}\leq n<\tau_{\zeta_n}^{(2, u)}\big)
\end{eqnarray*}

Based on (\ref{trebas}), (\ref{eqforcejump1}), \eqref{eqforcejump2-new} and (\ref{intermL}),  we only need  to estimate  $\p_{\eqref{mnx1234}}^{(4)}$.  By the many-to-one formula \eqref{manytoone},  we get that
\begin{eqnarray*}
\p_{\eqref{mnx1234}}^{(4)} &\leq& \e\Big[  \sum_{|u|=n}  1_{\{V(u)\leq \alpha_n-x,\,  \min_{\tau_{\zeta_n}^{(u)}\le j \le n}V(u_j) \geq \alpha_n-x-L ,\, \tau_{\zeta_n}^{(u)}\leq n<\tau_{\zeta_n}^{(2, u)}\}} \Big] 
\\
&\leq &\E\Big[  \ee^{S_n}  1_{\{S_n\leq \alpha_n-x,\,   \min_{\tau_{\zeta_n}\le j\le n}S_j\geq \alpha_n-x-L ,\, \tau_{\zeta_n} \leq n<\tau_{\zeta_n}^{(2)}\}} \Big] 
\\
&\leq &   \ee^{\alpha_n -x}\P\left( S_n\leq \alpha_n-x,\,   \min_{\tau_{\zeta_n}\le j\le n}S_j\geq \alpha_n-x-L,\, \tau_{\zeta_n} \leq n<\tau_{\zeta_n}^{(2)}\right)
\\
&\leq &  \ee^{\alpha_n -x}\sum_{k=1}^L\P\left( S_n-\alpha_n+x\in [-k,-k+1],\,    \min_{\tau_{\zeta_n}\le j\le n}S_j\geq \alpha_n-x-L,\, \tau_{\zeta_n} \leq n<\tau_{\zeta_n}^{(2)}\right)  \\
	&\leq & K \ee^{-x} L^2,
\end{eqnarray*}
where to obtain the last estimate, we have used several times the display (\ref{pasprov}) (with $y= -\alpha_n+x+L$, $a=  L-k$, $i\in [1,n]$ there). This  completes  the proof of Proposition \ref{upbound}.
\hfill$\Box$

\medskip

We end this section by a Lemma which will be used in Section \ref{S:mainestimate}:

\begin{Lemma}\label{L:n-t}  Assume (\ref{Hyp1})  and (\ref{Hyp2}).  Let  $\varepsilon>0$. For any $L>0$, there are some integers $T=T(\varepsilon, L)$ and $n_0=n_0(\varepsilon, L) >T $ such that for all $n\ge n_0$ and all $x>0$:  \begin{equation} \label{n-t} \p\Big( \exists u\in \mathbb{T}_n,\, V(u)\leq \alpha_n-x,\, \min_{\tau_{\zeta_n}^{(u)}\le j \le n}V(u_j)  \ge  \alpha_n-x -L,  \, \tau_{\zeta_n}^{(u)}\leq n- T, \,  \tau_{\zeta_n}^{(2, u)} > n \Big)\leq   \varepsilon\, \ee^{-x}.  \end{equation}
\end{Lemma}

{\noindent\it Proof of Lemma \ref{L:n-t} :}  Denote by $\p_{\eqref{n-t}}$ the probability term in \eqref{n-t} and write $y=\alpha_n -x -L$. Then by the many-to-one formula\begin{eqnarray*}
  \p_{\eqref{n-t}}  &\le& \e\left[  \sum_{\vert u \vert = n}  1_{\{  V(u)\leq y +L ,\, \min_{\tau_{\zeta_n}^{(u)}\le j \le n}V(u_j)  \ge y,  \, \tau_{\zeta_n}^{(u)}\leq n- T, \,  \tau_{\zeta_n}^{(2, u)} > n \}}  \right] \\
	&=& \E\left[ \ee^{S_n}   1_{\{  S_n \leq y+L,\, \min_{\tau_{\zeta_n} \le j \le n} S_j \ge y ,  \, \tau_{\zeta_n} \leq n- T, \, \tau_{\zeta_n}^{(2)} > n \}}  \right] \\
	&\le&  \sum_{i=1}^{n-T}  \sum_{k=1}^{\lfloor L\rfloor +1} \ee^{y+k} \P\Big( S_n  - y \in [k-1, k) ,\, \min_{\tau_{\zeta_n} \le j \le n} S_j \ge y,  \, \tau_{\zeta_n} =i, \, \tau_{\zeta_n}^{(2)} > n\Big) .
\end{eqnarray*}

\noindent Each probability term in the above double sum is less than, by \eqref{pasprov},  $K \P(\overline S_{n-i+1} \le k)   \, \ee^{-\alpha_n}  + n^{-(\alpha+\gamma-1)} \ell_3(n)  \le  K \P(\overline S_{n-i+1} \le L+1)   \, \ee^{-\alpha_n}  + n^{-(\alpha+\gamma-1)} \ell_3(n)  $, hence by taking the double sum over $i$ and $k$, $$   \p_{\eqref{n-t}} \le  K' \ee^{-x} \sum_{j=T}^n \P(\overline S_{j} \le L+1) + \ee^{-x} \, \ee^{\alpha_n}\, n^{- (\alpha+\gamma-2)} \ell_3(n). $$ 

Taking $T=T(\varepsilon, L)$ large enough such that $\sum_{j=T}^\infty  \P(\overline S_{j} \le L+1)  < {\varepsilon\over2K'}$ and $n_0$ large enough so that $\ee^{\alpha_n}\, n^{- (\alpha+\gamma-2)} \ell_3(n) \le {\varepsilon\over 2}$ for all $n \ge n_0$, we get \eqref{n-t}. \qed

%%%%%%%%%%%%%%%%
%%%%%%%%%%%%%%%% 
\section{Proof of  Proposition \ref{taildistrib}}\label{S:mainestimate}

At first we analyze  the trajectory of a particle which reaches  the minimum at time $n$.

 \medskip
  
Let $\varepsilon>0$ be small and $x >0$.  Let $L\equiv L(\varepsilon) \ge L_0 $ with $L_0$ is given by Lemma \ref{mauvaistraj}    be such that $K \ee^{- c_2 L }   < \varepsilon$.    Consider the event that there is some   $ u\in \mathbb{T}_n$ such that $ V(u)\leq \alpha_n-x$.   By   \eqref{eqforcejump1} and  \eqref{eqforcejump2-new}, with a cost at most $2 \varepsilon \ee^{-x}$, we may assume that $\tau_{\zeta_n}^{(u)}  \le n  $, which in view of \eqref{intermL2} yields that we  may furthermore assume   $ \min_{\tau_{\zeta_n}^{(u)}\le j \le n}V(u_j)  >   \alpha_n-x -L $  with an extra cost at most equal to~$\varepsilon \ee^{-x}$.  Finally by \eqref{n-t}, there exists some integer $T\equiv T(L , \varepsilon)$ such that we may assume $\tau_{\zeta_n}^{(u)} > n-T$ with an extra cost at most  equal to~$\varepsilon \ee^{-x}$.  Consequently   for all large $n \ge n_1(\varepsilon)$ and for all $x>0$,   \begin{eqnarray}
 &&\p\big( M_n\leq \alpha_n-x \big)   \nonumber 
 	\\ &= & \e\left[\frac{\sum_{|u|=n} 1_{ \{M_n=V(u) \leq \alpha_n-x,\,  \min_{\tau_{\zeta_n}^{(u)}\le j \le n}V(u_j)  \geq \alpha_n-x-L ,\,    \tau_{\zeta_n}^{(u)} \in [n-T, n]\} } }{\sum_{|u|=n} 1_{\{V(u)=M_n\}} }   \right] + O(\varepsilon) \, \ee^{-x}  \nonumber
	\\&=&\e_\q\left[ \ee^{V(w_n)} \frac{  1_{ \{M_n= V(w_n)\leq \alpha_n-x,\,  \min_{\tau_{\zeta_n}^{(w_n)} \le j \le n}V(w_j)  \geq \alpha_n-x-L ,\,     \tau_{\zeta_n}^{(w_n)} \in [n-T, n]\} } }{\sum_{|u|=n} 1_{\{V(u)=M_n\}} }  \right] + O(\varepsilon) \, \ee^{-x} , \label{bining1} 
\end{eqnarray}

\noindent where we have used the change of measure (cf. Proposition \ref{lyons}) for the last equality and  $O(\varepsilon)$ denotes, as usual,  some term bounded by a numerical constant times $\varepsilon$ (here by $5 \varepsilon$).

 \medskip

The next goal is to analyze the number of minima $\eta_n:= \sum_{|u|=n} 1_{\{V(u)=M_n\}}$ in  \eqref{bining1}. To this end,  we  consider  the  following event  \begin{equation}
\mathcal{E}_n(x):= \Big \{ \forall k  <  \tau_{\zeta_n}^{(w_n)},\, \forall v\in {\mathbb B}(w_k),\, \underset{u\geq v,\, |u|=n}{\min}V(u) >  \alpha_n-x\Big \},
\end{equation}
where  ${\mathbb B}(w_k) $, defined in \eqref{def-Omega}, denotes the set of brothers of $w_k$.  The following result  will be proved in Section \ref{S:aux}.

\begin{Proposition}  
\label{clust}
Under (\ref{Hyp1}), (\ref{Hyp2}) and (\ref{Hyp3}), for any $\varepsilon,\, L,\, T>0$ there exists $x_1>0$ such that for any $n\in \N$ large enough and $x\geq x_1$,
\begin{equation}
\label{eqclust} \q\left(   V(w_n)\leq \alpha_n-x ,\,    \min_{\tau_{\zeta_n}^{(w_n)} \le j \le n}V(w_j) \geq \alpha_n-x -L,\, \tau_{\zeta_n}^{(w_n)}\in [n-T, n],\,   (\mathcal{E}_n(x))^c\right)\leq \varepsilon\,     \ee^{-\alpha_n}  .
\end{equation}
Consequently for all $x\ge x_1$ and all large $n$, \begin{equation}
\label{eqclust2}\mathbb E_\q\left[ \ee^{V(w_n)} 1_{ \{   V(w_n)\leq \alpha_n-x ,\,   \min_{\tau_{\zeta_n}^{(w_n)} \le j \le n}V(w_j) \geq \alpha_n-x -L,\, \tau_{\zeta_n}^{(w_n)}\in [n-T, n],\,   (\mathcal{E}_n(x))^c \}} \right]\leq \varepsilon\,  \ee^{-x} .
\end{equation}
\end{Proposition}

\medskip

 For any $u \in \T$ and    $j\ge0$,  we define \begin{equation}\label{eta} M_j(u):= \min_{ |v|= |u|+j, \, v  \in \T^{(u)}} (V(v)- V(u)) , \qquad  \eta_j(u):= \sum_{ |v|= |u|+j, \, v \in \T^{(u)}} 1_{\{ V(v)- V(u) = M_j(u)\}} ,
\end{equation}

\noindent  with $ M_j(\varnothing) \equiv M_j$  and $\eta_j(\varnothing)\equiv \eta_j$.   In the case that the subtree $\T^{(u)}$ does not survive up to $j$-th generation,   by definition   $M_j(u)=\infty$ and   $\eta_j(u)=0$ [which is in agreement with the convention that $\sum_\emptyset\equiv 0$].

\medskip

On the event $ \mathcal{E}_n(x) \cap \{M_n\leq \alpha_n-x\}$,   any particle located at the minimum stays on the spine at least up to  the generation $ \tau_{\zeta_n}^{(w_n)} -1 $. Therefore  on $\mathcal{E}_n(x) \cap\{\tau_{\zeta_n}^{(w_n)}=k\}$,   for any $|u|=n$ satisfying  that $V(u)=  M_n$,  there is some $v$ with ${\buildrel \leftarrow \over v}=  w_{k-1}$  such  that  $u \in \T^{(v)}$ (either $v=w_k $ or   $v \in {\mathbb B}(w_k)$).   We obtain that on  $  \mathcal{E}_n(x) \cap\{\tau_{\zeta_n}^{(w_n)}=k\}$  with   $k \le n$,  \begin{eqnarray}  \eta_n \equiv \sum_{|u|=n}  1_{\{V(u)=M_n\}}  &  = &     \sum_{{\buildrel \leftarrow \over v}=  w_{k-1}}   \sum_{|u|=n , u \in \T^{(v)}}  1_{\{V(u)=M_n\}}   \nonumber
	\\ &=&   \sum_{{\buildrel \leftarrow \over v}=  w_{k-1}}   \eta_{n-k}(v) \,  1_{\{ M_{n-k}(v)= M_n - V(v)\}}  .   \label{eta2}  \end{eqnarray}

 In view of  \eqref{bining1} and \eqref{eqclust2}, we deduce from \eqref{eta2}  that   for any   $ x\geq x_1$,  for  $n$ large enough,
\begin{eqnarray}
 && \mathbb P(M_n\leq \alpha_n-x) \nonumber \\
 	&&=\mathbb E_\q\Big[ {\ee^{V(w_n)} \over \eta_n} 1_{\{M_n= V(w_n)\leq \alpha_n-x,\,  \min_{\tau_{\zeta_n}^{(w_n)} \le j \le n}V(w_j)  \geq \alpha_n-x-L,\,  \tau_{\zeta_n}^{(w_n)}\in [ n-T,n]\}} , \,  \mathcal{E}_n(x)  \Big]  + O(\varepsilon) \ee^{-x}  
	\nonumber \\ &&=  \sum_{k= n-T}^n \mathbb E_\q\Big[ \ee^{V(w_n)}\frac{1_{\{M_n= V(w_n)\leq \alpha_n-x,\,   \min_{k  \le j \le n}V(w_j)  \geq \alpha_n-x-L,\,  \tau_{\zeta_n}^{(w_n)} = k \}} }{  \sum_{{\buildrel \leftarrow \over v}=  w_{k -1}}   \eta_{n-k }(v) \,  1_{\{ M_{n-k}(v)= V(w_n) - V(v)\}} },\,  \mathcal{E}_n(x)  \Big]  + O(\varepsilon) \ee^{-x}  \nonumber \\
	&&=  \sum_{k= n-T}^n  \mathbb E_\q\Big[ A_{\eqref{mn-x1}} (k)  \Big]  + O(\varepsilon) \ee^{-x}  ,  \label{mn-x1}
\end{eqnarray}

\noindent where   \begin{eqnarray*}
A_{\eqref{mn-x1}} (k)&:=& \ee^{V(w_n)}\frac{1_{\{M_n= V(w_n)\leq \alpha_n-x,\,   \min_{k  \le j \le n}V(w_j)  \geq \alpha_n-x-L,\,  \tau_{\zeta_n}^{(w_n)} = k \}} }{  \sum_{{\buildrel \leftarrow \over v}=  w_{k -1}}   \eta_{n-k }(v) \,  1_{\{ M_{n-k}(v)= V(w_n) - V(v)\}} }    , \qquad n-T \le k \le n, 
\end{eqnarray*}

\noindent and  the last equality in \eqref{mn-x1}  still holds thanks to   \eqref{eqclust2}.  Obviously  the following upper bound   holds:   \begin{equation} \label{A-upp}
A_{\eqref{mn-x1}} (k) \le \ee^{V(w_n)}\, 1_{\{ V(w_n)\leq \alpha_n-x,  \min_{k  \le j \le n}V(w_j)  \geq \alpha_n-x-L, \tau_{\zeta_n}^{(w_n)} = k\}} =:  B_{\eqref{A-upp}}(k)
\end{equation}

\noindent 
 Moreover, under $\q$, $(V(w_j), j\ge0)$ is distributed as the random walk  $(S_j, j\ge0)$ under $\P$.  Then  \begin{equation} \label{expect-B}
\e_\q \left[B_{\eqref{A-upp}}(k)  \right] \le  \ee^{\alpha_n -x} \P\Big(  S_n\le \alpha_n -x,  \min_{k\le j \le n} S_j  \ge \alpha_n -x-L , \tau_{\zeta_n} = k\Big) \le K\, \ee^{-x},
\end{equation}

\noindent by using \eqref{pasprov-sum}.

In view of the hypothesis \eqref{decorr} which  only holds for those functions with compact support, we need to  truncate  $| \Delta V(w_k)- \Delta V(v)| $ uniformly on $v \in {\mathbb B}(  w_{k }) $. This is possible thanks to    the following~Claim:

\begin{Claim}\label{Claim-compact} There exists some $\lambda_0=\lambda_0(\varepsilon, L, T)>0$ such that for all $\lambda\ge \lambda_0$, if we define  the event  $ \Upsilon_k(\lambda,  n, T) $ by   $$\Upsilon_k(\lambda,  n, T) ^c  :=   \Big\{   \exists v \in {\mathbb B}(w_k):  | \Delta V(w_k) - \Delta V(v) | > \lambda,  M_{n-k}(v) \le   V(w_n) - V(v) \Big \}  ,  \qquad n-T \le k \le n, $$  then $$  \sum_{k= n-T}^n\mathbb E_\q\Big[ B_{\eqref{A-upp}}(k) , \Upsilon_k(\lambda,  n, T)  ^c  \Big]   \le  O(\varepsilon)\, \ee^{-x}, $$
in particular, $$  \sum_{k= n-T}^n\mathbb E_\q\Big[ A_{\eqref{mn-x1}} (k)  , \Upsilon_k(\lambda,  n, T)  ^c  \Big]   \le  O(\varepsilon)\, \ee^{-x}. $$
\end{Claim}

 \medskip
 Observe that  on  $\Upsilon_k(\lambda, n, T)$,  for any $v \in {\mathbb B}(w_k)$ satisfying that $ | \Delta V(w_k) - \Delta V(v) |  >\lambda$, we have that $M_{n-k}(v) >  V(w_n) - V(v)$, hence   the subtree $\T^{(v)}$  contains a possible (global) minimum only if $v \in {\mathbb B}_\lambda(w_k)$, where
 \begin{equation}
\label{Omega-lambda} {\mathbb B}_\lambda(w_k):= \Big\{ v \in {\mathbb B}(w_k): |  \Delta V(w_{k}) - \Delta V(v) | \le  \lambda \Big\}, \qquad n-T \le k \le n. 
\end{equation}
It follows that on  $\Upsilon_k(\lambda, n, T)$, 
	\begin{eqnarray} &&  \sum_{{\buildrel \leftarrow \over v}=  w_{k -1}}   \eta_{n-k }(v) \,  1_{\{ M_{n-k}(v)= V(w_n) - V(v)\}}   \nonumber\\
	&&=   \eta_{n-k}(w_k) +  \sum_{ v \in {\mathbb B}_\lambda(w_k)}   \eta_{n-k }(v) \,  1_{\{ M_{n-k}(v)= V(w_n) -V(w_k) + \Delta V(w_k)- \Delta V(v)\}}  \nonumber \\
	&&=: I_{\eqref{eta-decom}} , 
	\label{eta-decom} \end{eqnarray}

\noindent whereas    \begin{eqnarray}
1_{\{M_n= V(w_n) \} } &=& 1_{ \{ 	M_{n-k}(w_k)= V(w_n)-V(w_k)  \}  }      \prod_{  v \in {\mathbb B}_\lambda(  w_k )} 1_{\{    M_{n-k}(v) \ge V(w_n) -V(w_k) + \Delta V(w_k)- \Delta V(v)\}}   \nonumber\\
	& =: & J_{\eqref{mn-decom}}.  \label{mn-decom}
\end{eqnarray}

 Therefore on $\Upsilon_k(\lambda, n, T)$,  \begin{equation} \label{A-lam}
 A_{\eqref{mn-x1}} (k) =  B_{\eqref{A-upp}}(k)  { J_{\eqref{mn-decom}}  \over I_{\eqref{eta-decom}}}, 
\end{equation}

\noindent which  is the key   to truncate the point measure  $\sum_{v\in {\mathbb B} (w_k)} \delta_{\{ V(w_k) - V(v)\}} $  to   $\sum_{v\in {\mathbb B}_\lambda (w_k)} \delta_{\{ V(w_k) - V(v)\}} $.

 \medskip
 
 Before using  \eqref{A-lam},  we give the proof of Claim \ref{Claim-compact}.

  \medskip
 {\noindent \it Proof of Claim \ref{Claim-compact}.} By \eqref{tight+},   there exists some $\lambda_1>0$ and some $z_0\in \r$  such that for all  $\lambda\ge \lambda_1$ and $z \le z_0$,  \begin{equation} \label{tight+1}
   \q \Big( \cup_{ v \in {\mathbb B}(w_k)} \{ \Delta V(w_k) - \Delta V(v) \ge \lambda\} \,   \big | \, \Delta V(w_k)= z \Big)   \le {\varepsilon\over T},
		\end{equation}

\noindent where we remark that the probability term in \eqref{tight+1} does not depend on $k$.   Observe  that  for all large $n$ [such that $ - \zeta_n \le z_0$], on $\{ \Delta V(w_k)\leq -\zeta_n \}$,  $ \q \Big( \cup_{ v \in {\mathbb B}(w_k)} \{ \Delta V(w_k) - \Delta V(v) \ge \lambda\} \,   \big | \, B_{\eqref{A-upp}}(k)  \Big) =\q \Big( \cup_{ v \in {\mathbb B}(w_k)} \{ \Delta V(w_k) - \Delta V(v) \ge \lambda\} \,   \big | \, \Delta V(w_k)  \Big) \le {\varepsilon\over T} $ by  \eqref{tight+1}. It follows that   \begin{equation}\label{tight+2}
 \mathbb E_\q \left[ B_{\eqref{A-upp}}(k),  \cup_{ v \in {\mathbb B}(w_k)} \{ \Delta V(w_k) - \Delta V(v) \ge \lambda\} \right]  \nonumber\\  
 \le   {\varepsilon\over T} \mathbb E \left[  B_{\eqref{A-upp}}(k) \right]    \le  K{\varepsilon\over T}  \ee^{-x}.  \end{equation}

 Write $\#{\mathbb B}(w_k) = \sum_{v \in {\mathbb B}(w_k)} 1$.   By \eqref{tightness} we may choose a large constant $\lambda_2$  and some $z_0< 0$ such that for all $z \le z_0$,  $$ 
\q\Big( \# {\mathbb B}(w_k) > \lambda_2 \, \big |\,    \Delta V(w_k) = z \Big) \le  { \varepsilon\over T}, \qquad \forall \, k \ge 1.$$
We notice that the above probability  does not depend on $k$.

\medskip

Now, we treat  the case $\Delta V(w_k) - \Delta V(v)  < - \lambda$: At first,  \begin{eqnarray}
 &&\mathbb E_\q\Big[ B_{\eqref{A-upp}}(k) , V(w_n) - V(w_k) \ge \lambda_3  \Big]   \nonumber \\
  &\le &     \ee^{\alpha_n -x}\, \P \Big( S_n\le \alpha_n -x,  \min_{k\le j \le n} S_j\ge \alpha_n -x-L , \tau_{\zeta_n} = k,  S_n - S_k   >  \lambda_3  \Big) \nonumber 
  \\ &\le& { \varepsilon\over T}  \, \ee^{-x}, \label{tight-1}
\end{eqnarray}

\noindent by applying \eqref{tight2}  and by choosing    a constant $\lambda_3=\lambda_3( \varepsilon, T, L) $  large enough.  For those  $v \in {\mathbb B}(w_k)$ such that the event $\Upsilon_k(\lambda, n, T)^c$ holds, if furthermore    $\Delta V(w_k) - \Delta V(v)  < - \lambda $ and $V(w_n) - V(w_k) < \lambda_3$,  then $M_{n-k} (v)  < - \lambda + \lambda_3$ which holds with a probability bounded from above by $\ee^{-\lambda+\lambda_3}$ (see Lemma~\ref{L:expmin}).

\medskip

Recall that under $\q$,    $(\Delta V(w_j), \sum_{v \in {\mathbb B}(w_j)} \delta_{\{ \Delta V(v)\}}, j\ge 1)$ is a sequence of i.i.d. random variables, whereas conditioning on ${\cal G}$,   for   $v \in {\mathbb B}(w_k)$, $ (\eta_{n-k}(v), M_{n-k}(v))$ are independent and are distributed as $(\eta_{n-k}, M_{n- k}) $ under $\mathbb P$.  It follows from \eqref{tight-1} that  
\begin{eqnarray}
&&\mathbb E_\q\Big[ B_{\eqref{A-upp}}(k) ,  \exists v \in {\mathbb B}(w_k):    \Delta V(w_k) - \Delta V(v)   < -  \lambda,  M_{n-k}(v) \le   V(w_n) - V(v)    \Big]   \nonumber  \\
	&&\le { \varepsilon\over T}  \, \ee^{-x} + \mathbb E_\q\Big[ B_{\eqref{A-upp}}(k)  \times  \big(  1_{\{ \# {\mathbb B}(w_k) > \lambda_2\}} +    \lambda_2\times    \ee^{- \lambda + \lambda_3} \big) \Big]  \nonumber \\
	&&\le  { \varepsilon\over T}  \, \ee^{-x}  + \lambda_2     \ee^{- \lambda + \lambda_3} \mathbb E_\q\Big[ B_{\eqref{A-upp}}(k) \Big]+\mathbb E_\q\Big[ B_{\eqref{A-upp}}(k)  1_{\{ \# {\mathbb B}(w_k) > \lambda_2\}}  \Big]  =: C_{\eqref{tight-fin}}. \label{tight-fin}  
\end{eqnarray}
 
  Since  $\Delta V(w_k) \le - \zeta_n \le z_0, $ we have $\q\Big( \# {\mathbb B}(w_k) > \lambda_2 \, \big |\,  {\cal G}      \Big)= \q\Big( \# {\mathbb B}(w_k) > \lambda_2 \, \big |\,    \Delta V(w_k)   \Big) \le {\varepsilon\over T}$. Then we deduce  from \eqref{expect-B}   that  $$ C_{\eqref{tight-fin} } \le { \varepsilon\over T}  \, \ee^{-x}  + \lambda_2     \ee^{- \lambda + \lambda_3}   K  \ee^{-x}+ {\varepsilon\over T} K \ee^{-x} = O(\varepsilon) \ee^{-x}, $$

  \noindent for all $\lambda $ large enough [$\lambda_2$ and $\lambda_3$ being fixed].  This and \eqref{tight+2} yield  Claim \ref{Claim-compact}. 
 \qed
 
 	\medskip

	Based on Claim \ref{Claim-compact} and \eqref{A-lam}, we obtain that  for all $\lambda\ge \lambda_0(\varepsilon, L, T) $, $n$ large enough and $x\ge x_1$:
\begin{eqnarray} \label{new189} 
 \mathbb P(M_n\leq \alpha_n-x) & = & \sum_{k= n-T}^n  \mathbb E_\q\Big[ A_{\eqref{mn-x1}} (k)  , \, \Upsilon_k(\lambda, n, T) \Big]  + O(\varepsilon) \ee^{-x}  \nonumber \\
 &=&  \sum_{k= n-T}^n  \mathbb E_\q\Big[ B_{\eqref{A-upp}}(k)  { J_{\eqref{mn-decom}}  \over I_{\eqref{eta-decom}}} \, \Upsilon_k(\lambda, n, T) \Big]  + O(\varepsilon) \ee^{-x}  \nonumber \\
 	&=& \sum_{k= n-T}^n  \mathbb E_\q\Big[ B_{\eqref{A-upp}}(k)  { J_{\eqref{mn-decom}}  \over I_{\eqref{eta-decom}}}   \Big]  + O(\varepsilon) \ee^{-x} , 
\end{eqnarray}
 
\noindent   by using again Claim \ref{Claim-compact}.  	Write again for brevity $$ y\equiv y(n, x, L):= \alpha_n -x -L.$$

 Observe that   $\{ \tau_{\zeta_n}^{(w_n)} = k\}= \{  \tau_{\zeta_n}^{(w_k)}= k\}$.  In view of  \eqref{eta-decom} and \eqref{mn-decom}, we deduce from the Markov property at $k$  that    \begin{eqnarray}\label{bmnx}
   && \mathbb E_\q\Big[ B_{\eqref{A-upp}}(k)  { J_{\eqref{mn-decom}}  \over I_{\eqref{eta-decom}}}   \Big] \nonumber 
   \\&=&    \ee^{\alpha_n -x} \,\mathbb E_\q\left[  1_{ \{ \tau_{\zeta_n} ^{(w_k)}=k,   V(w_k) \ge  y \}}    F^{(L)}_{n-k} \Big(  V(w_k)- y  , \sum_{v \in {\mathbb B}_\lambda(w_k)} \delta_{\{ \Delta V(w_{k})- \Delta V(v)\}} \Big)  \right] , \end{eqnarray}

  \noindent where for any $j \ge0$, $F^{(L)}_j: \r_+ \times   {\cal M} \to  [0, 1] $ is the measurable function defined  as follows: for any $\theta \in  {\cal M}$, say $\theta= \sum_{i=1}^l \delta_{\{y_i\}}$ with some $l\ge1$ and $|y_i|\le \lambda$, \begin{eqnarray} &&  \label{def-fjl}   F_j^{( L)}(s,   \theta)  \nonumber 
  	\\ &: =&  \ee^{s-L} \,\mathbb E_\q \Big[ \ee^{V(w_j)} {1_{\{ M_j= V(w_j )\}} \,  1_{\{V(w_j) \le L-s , \underline V(w_j) \ge  -s   \}} \over \eta_j  + \sum_{i=1}^l \eta_j^{(i)} 1_{\{ M_j^{(i)} =V(w_j) +  y_i  \}}}    \prod_{1\le i \le l, \, \eta_j^{(i)}\neq 0} 1_{\{ M_j^{(i)} \ge V(w_j) + y_i  \}} \Big] ,  \end{eqnarray}

  \noindent where $ \underline V(w_j):= \min_{0\le i \le j} V(w_i)$ and  under $\q$, $( \eta_j^{(i)}, M_j^{(i)}, j\ge0)_{i\ge1}$ is   i.i.d., independent of everything else and distributed as $(\eta_j , M_j, j\ge0)$ under $\p$.  We mention that if $\eta_j^{(i)}=0$, then $M_j^{(i)}=\infty$ by definition. Obviously, the above expectation under $\q$ does not depend on the order of $\{y_i\}$ in $\theta$.    Recall that  under $\q$, $(\Delta V(w_k), \sum_{v \in {\mathbb B}(w_k)} \delta_{\{ \Delta V(v)\}})_{k\ge1}$ are i.i.d., and $(V(w_j), j\ge 0)$ is distributed as $(S_j, j\ge 0)$ under $\P$.  If we define     \begin{equation}\label{def-gjl} G_j^{(\lambda, L)}(s, z):=\mathbb E_\q \Big[ F^{(L)}_j \Big(s,   \sum_{v \in {\mathbb B}_\lambda(w_1)} \delta_{\{   V(w_{1})-   V(v)\}} \Big)  \big| \,   V(w_1)= z\Big], \qquad j\ge 0, \, s, z \in \r, \end{equation}

  \noindent then      \begin{equation}\label{bmnx2}
\mathbb E_\q\Big[ B_{\eqref{A-upp}}(k)  { J_{\eqref{mn-decom}}  \over I_{\eqref{eta-decom}}}    \Big]  =  \ee^{\alpha_n -x}\,    \E\left[  1_{ \{ \tau_{\zeta_n} =k,   S_k\ge  y  \}}   \, G^{(\lambda, L)}_{n-k} \Big(   S_k- y   ,   X_k  \Big)  \right] ,\end{equation}

  \noindent where as before, $X_k= S_k- S_{k-1}$.   For any $j\ge 0$,  notice that $\eta_j \ge 1_{\{V(w_j)= M_j\}}$, hence $ F_j^{(\lambda, L)}(s, z, \theta) \le   \ee^{s-L} \,\mathbb E_\q \left[ \ee^{V(w_j)}     1_{\{V(w_j) \le L-s , \underline V(w_j) \ge  -s \}}   \right] \le  \q \big( V(w_j) \le L-s \big)= \P(S_j \le L-s).$   It follows that   for any $s\ge 0$ and $z \in \r$,    \begin{equation}\label{renew-cond1}  G^{(\lambda, L)}_j(s, z) \le \P(S_j \le L-s)   \le j \, \P( X \le {L- s \over j} ) \le K_{j, L} (1+s)^{- \alpha}, 
  \end{equation}
 with some constant $K_{j, L}>0$. 
  
 \medskip
 
 Recalling (\ref{decorr}), let $\Xi= \sum_{i=1}^{\nu^*} \delta_{ \{ y^*_i\}}$, $y^*_i \in \r \cup\{-\infty\}$  be a point process independent of everything else whose law is defined as the limiting law of $\sum_{w\in \mathbb{B}(w_1)}\delta_{\{V(w_1)-V(v)\}}$ conditioned on $\{V(w_1)=z\}$ as $z \to -\infty$.  We claim that    as $z \to -\infty$,  for any $s\ge0$,  \begin{equation} \label{renew-cond2} G^{(\lambda, L)}_j(s, z) \mbox{ converges to }  G^{(\lambda, L)}_j(s),  \end{equation}  

\noindent where    \begin{eqnarray*}
&& G^{(\lambda, L)}_j(s) \\ &:=& \ee^{s- L} \,\mathbb E_\q \Big [ \ee^{V(w_j)} {1_{\{ M_j= V(w_j) \}} \,  1_{\{ V(w_j)\le  L-s ,\,  \underline V(w_j) \ge - s   \}} \over \eta_j  + \sum_{1\le i \le \nu^*: |y_i^*|\le \lambda }  \eta_j^{(i)} 1_{\{ M_j^{(i)} =V(w_j)  +  y^*_i \}}}      \prod_{1\le i \le \nu^*:  \eta_j^{(i)}\neq0, |y^*_i| \le \lambda }1_{\{ M_j^{(i)} \ge V(w_j) + y^*_i   \}} \Big]  , 
\end{eqnarray*}

\noindent with the usual convention $\prod_\emptyset:= 1$.  Moreover,  \begin{eqnarray*}
\int_0^\infty G_j^{(\lambda, L)}(s) ds &=& \mathbb E_\q \Big[ \big( \ee^{- V(w_j)}  - \ee^{- L- \underline V(w_ j)}\big) {1_{\{ M_j= V(w_j) \}} \,  1_{\{V(w_j)\le  L \}  }  \over \eta_j  + \sum_{1\le i \le \nu^*:  |y^*_i| \le \lambda}  \eta_j^{(i)} 1_{\{ M_j^{(i)} =V(w_j) + y^*_i \}}}  \nonumber 
	\\ && \qquad\qquad \prod_{1\le i \le \nu^*: \eta_j^{(i)}\neq0, |y^*_i| \le \lambda}1_{\{ M_j^{(i)} \ge V(w_j) + y^*_i ,   | y^*_i | \le   \lambda \}} \Big] .
\end{eqnarray*}

   Let us postpone for the moment the proof of \eqref{renew-cond2}.   
    By  assembling \eqref{new189} and  \eqref{bmnx}, we get that for any $x\ge x_1$, $L\ge L_1$ and large $T$, for all large $n\ge n_0(\varepsilon, L, T)$,  \begin{equation}
\mathbb P(M_n\leq \alpha_n-x)  =    \ee^{\alpha_n -x }\, \sum_{j=0}^{T-1} \E\left[  1_{ \{ \tau_{\zeta_n} =n-j, S_{n-j} \ge y\}} \,  G_j^{(\lambda,L)}(  S_{n-j} - y, \ X_{n-j})  \right] + O(\varepsilon) \, \ee^{-x}    .  \label{bining2} 
\end{equation}

By means of \eqref{renew-cond1} and \eqref{renew-cond2}, we can    apply  Lemma \ref{renewal} to  $G^{(\lambda, L)}_j(s, z)$, for any fixed $0\le j < T$.   This  gives that    for  any    $L\ge L_1$ and large $T$, for all large $n\ge n_1(\varepsilon, L, T)$  and $x\in [ x_1,\, \frac{n}{\log n}]$,    \begin{equation}
\label{sump} \left |\mathbb P\left(M_n\leq \alpha_n-x\right)- \mathtt{m}^{-(\alpha+1)}\ee^{-x}  \sum_{j=0}^{T-1} \int_{0}^\infty G^{(\lambda, L)}_{j}(s)ds  \right|\leq  O(\varepsilon)\,  \ee^{-x}.
\end{equation}

On the other hand,  we deduce from the  bounded convergence theorem (when $\lambda$ tends to $\infty$) and monotone convergence  theorem (when $T$ and $L$ tend to $\infty$) that  
\begin{eqnarray}
\label{defC*}  && \underset{L\to\infty}{\lim}\, \underset{T\to \infty}{\lim}\,  \lim_{\lambda\to\infty}   \mathtt{m}^{-(\alpha+1)}  \sum_{j=0}^{T-1} \int_{0}^\infty  G^{(\lambda, L)}_j(u)du  \nonumber
	\\ && = \mathtt{m}^{-(\alpha+1)} \sum_{j=0}^\infty\mathbb E_\q \left[  \ee^{- V(w_j)}  {1_{\{ M_j= V(w_j) \}}    \over \eta_j  + \sum_{i=1}^{\nu^*}  \eta_j^{(i)} 1_{\{ M_j^{(i)} =V(w_j) +   y^*_i \}}} \prod_{i=1}^{\nu^*}1_{\{ M_j^{(i)} \ge V(w_j) + y^*_i   \}} \right]  \nonumber
		\\&&   =:c_*,
\end{eqnarray}
[in the product $ \prod_{i=1}^{\nu^*}$, if $\eta_j^{(i)}=0$ then $M_j^{(i)}=\infty$ by definition], moreover by Proposition \ref{upbound} we know that $c_*$ is a finite constant. By combining (\ref{sump}) and (\ref{defC*}) we get  Proposition \ref{taildistrib}. 

\medskip

It remains to check \eqref{renew-cond2}.     If we denote by ${\bf  e}$ an independent standard exponential variable, then we may rewrite \eqref{def-fjl} as \begin{eqnarray*}  F_j^{(L)}(s,  \theta)  
	 & =&  \ee^{s-L} \mathbb E_\q \Big[ \ee^{V(w_j)}  1_{\{ M_j= V(w_j )\}} \,  1_{\{V(w_j) \le L-s , \underline V(w_j) \ge  -s   \}}  \ee^{-  {\bf e} (\eta_j-1) }  \\
	&& \qquad \qquad \times   \prod_{1\le i \le l, \eta_j^{(i)}\neq 0}  1_{\{ M_j^{(i)} \ge V(w_j) + y_i      \}} \ee^{- {\bf e}  \,    \eta_j^{(i)} 1_{\{ M_j^{(i)} =V(w_j) +  y_i \}}}   \Big] \\
	&=&  \ee^{s-L} \mathbb E_\q \left[ \ee^{V(w_j)}  1_{\{ M_j= V(w_j )\}} \,  1_{\{V(w_j) \le L-s , \underline V(w_j) \ge  -s   \}}  \ee^{-  {\bf e} (\eta_j-1) -      \sum_{i=1}^l    h^{(\lambda)}_{ j, {\bf e},  V(w_j) } ( y_i)  } \right] ,  \end{eqnarray*}

\noindent where    we have used the fact that $( \eta_j^{(i)}, M_j^{(i)}, j\ge0)_{i\ge1}$ is   i.i.d., independent of everything else and distributed as $(\eta_j , M_j, j\ge0)$ under $\mathbb P$, and  for any $ a>0, b \in \r $ and $x \in \r$ [remark that $|y_i | \le \lambda$ when $\theta= \sum_{v\in \mathbb{B}_\lambda(w_1)}\delta_{\{V(w_1)-V(v)\}}$], $$ h^{(\lambda)}_{j ,a, b }(x):= -  1_{\{ |x| \le \lambda\}}\,  \log  \mathbb E\Big[ 1_{\{ \eta_j=0\}} + 1_{\{ \eta_j \neq  0, M_j \ge b+x \}} \ee^{- a \eta_j 1_{\{ M_j= b+ x\}}} \Big]   .$$ 

 Then \eqref{renew-cond2} follows from the assumption \eqref{decorr} and an application of dominated convergence theorem.  The proof of Proposition \ref{taildistrib} is now completed. 
\hfill$\Box$

\section{Proof of Proposition  \ref{clust}} \label{S:aux}

Fix $0<\varrho<\min(\frac{\alpha-1}{2}, \frac{1}{12})$.  Recall \eqref{def-Omega} that  ${\mathbb B}(u)$ is the set of   brothers of $u$ for any $u \in \T\backslash\{\varnothing\}$.    Let $B>0$ be a large constant and $J$ be a large integer.  Recall  \eqref{drop-u1} and \eqref{drop-u2}. 

Let us say that  $u \in \T_n$ is  a   good  vertex  if for any $x\geq 0$, 
\begin{equation}
\label{defGood}
\tau_{\zeta_n}^{(2,u)}>n \ge \tau_{\zeta_n}^{(u)} > J,\quad \text{and }\quad \sum_{v\in {\mathbb B}(u_k)}\ee^{-(V(v)+x)}\leq \left\{ \begin{array}{ll} \ee^{B-x},\quad &\text{if }  1 \le k \le J,
\\
\ee^{-k^\varrho},\quad &\text{if }  J < k <   \tau_{\zeta_n}^{(u)}. 
\end{array}\right.
\end{equation}

% Define  \begin{equation} \label{defGood} {\mathbb G}_{n, J, B}:= \Big\{ \tau_{\zeta_n}^{(2, w_n)}>n \ge \tau_{\zeta_n}^{(w_n)}> J  \Big\} \cap \Big\{ \sum_{v\in {\mathbb B}(u_k)}\ee^{-(V(v)+x)}\leq  1_{\{1\le k \le J\}} \ee^{B-x} + 1_{\{ J< k <  \tau_{\zeta_n}^{(w_n)}\}} \ee^{-k^{\varrho}}  \Big\}. \end{equation}

The condition $\{n \ge \tau_{\zeta_n}^{(u)}> J \}$ will be automatically satisfied in the event that we are interested in.   
Roughly saying,  when $w_n$ is good,  the contribution from the particles in ${\mathbb B}(w_k)$, for all $k <  \tau_{\zeta_n}^{(w_n)}$, is not too large.  The following lemma estimates the case when $w_n$  is not good:

%Let us   say that $|u|=n$ is a good vertex if \begin{equation} \label{defGood} \tau_{\zeta_n}^{(2)}(u)>n,\quad \text{and }\quad \sum_{v\in {\mathbb B}(u_k)}\ee^{-(V(v)+x)}\leq \left\{ \begin{array}{ll} \ee^{B-x},\quad &\text{if } k< \log B, \\ \ee^{-k^\frac{1}{12}},\quad &\text{if } k\in [\log B,\, \tau_{\zeta_n}(u)]. \end{array}\right. \end{equation}

\begin{Lemma}
\label{L:App1}
  Under (\ref{Hyp1}), (\ref{Hyp2}) and (\ref{Hyp3}), for any $L, T,\, \varepsilon >0$,  there exists $J(L,T,\varepsilon)$ such that for all $J\ge J(L,T,\varepsilon)$ , there exists $B(J,L,T,\varepsilon)>0$ such that for all $B\ge B(J,L,T,\varepsilon)$,  for any $n$ large enough and  $ x\geq 0$,
\begin{equation}
\label{eqL:App1}
\q\Big(V(w_n)\leq \alpha_n-x ,\, \min_{\tau_{\zeta_n}^{(w_n)} \le j \le n}V(w_j)\geq \alpha_n-x -L,\, \tau_{\zeta_n}^{(w_n)}\in  [n-T,n],\,   \mbox{ $w_n$  not good }  \Big)\leq \varepsilon    \, \ee^{-\alpha_n} .
\end{equation}
\end{Lemma}

\medskip
By admitting Lemma \ref{L:App1} for the moment, we give now the proof of Proposition \ref{clust}: 
\medskip

\noindent{\bf Proof of Proposition \ref{clust}.}  For brevity we use the following notation: \begin{equation} \label{def-fn}
 F_n \equiv F_{n, T, x} : = \Big\{ V(w_n)\leq \alpha_n-x ,\, \min_{\tau_{\zeta_n}^{(w_n)} \le i \le n}V(w_j) \geq \alpha_n-x -L,\, \tau_{\zeta_n}^{(w_n)}\in  [n-T,n]\Big\} 
\end{equation}

 By \eqref{eqL:App1}, it remains to estimate the following probability: 
\begin{eqnarray}    \mathbb Q_{\eqref{goodwn}} \label{goodwn}  &:= &  \q\Big( F_n, \, w_n \text{ good},\,   (\mathcal{E}_n(x))^c\Big)  \nonumber
	\\ & \le &\mathbb E_\q\Big[   1_{\{  F_n,\, w_n \text{ good}\}}\sum_{j=1}^{\tau_{\zeta_n}^{(w_n)}-1}  \sum_{v\in {\mathbb B}(w_j)}     1_{\{ \underset{u\geq v,\, |u|=n}{\min} V(u)\leq  \alpha_n - x\}} \Big] \nonumber
\\
& \le &  \sum_{t=n-T}^n  \sum_{j=1}^{t-1}\mathbb E_\q\Big[ 1_{\{  F_n , \, t=\tau_{\zeta_n}({w_n})  ,\,  w_n \text{ good}  \}}   \sum_{v\in {\mathbb B}(w_j)}    1_{\{ \underset{u\geq v,\, |u|=n}{\min} V(u)\leq  \alpha_n - x\}}    \Big] .  
\end{eqnarray}

By the spinal decomposition (Proposition \ref{lyons} (iii)),  for any $t\in [n-T,n]$, $j\in [1,t-1]$ and $v\in {\mathbb B}(w_j)$,  conditionally on ${\cal G}$ and on $\{ V(v)= b\}$, we have 
 $$  
 \q\Big(     \underset{u\geq v,\, |u|=n}{\min} V(u)\leq \alpha_n - x   \big | \, {\cal G}\Big)=      \mathbb P\big( M_{n- j }\leq \alpha_n-x- b  \big) .  
 $$

If $j \le {2n\over3}$, we apply Proposition \ref{upbound} to get that $  \mathbb  P\big( M_{n- j }\leq \alpha_n-x- b  \big) \le K\, \ee^{- ( b+x +\alpha_{n-j}- \alpha_n)}$, whereas if ${2n \over 3} < j \le t$, we apply Lemma \ref{L:expmin} (which holds obviously for all $x\in \r$) and get that $  \mathbb P\big( M_{n- j }\leq \alpha_n-x- b  \big) \le \ee^{- (b+x- \alpha_n)}$.  Taking into  account the fact that $w_n$ is good, we obtain  \begin{eqnarray*}
\E_\q\Big[  \sum_{v\in {\mathbb B}(w_j)}  1_{\{ \underset{u\geq v,\, |u|=n}{\min} V(u)\leq \alpha_n - x\}}   \big | \, {\cal G} \Big]  \leq \left\{ \begin{array}{ll}  2 K\, \ee^{ \alpha_n-\alpha_{n- j } } \ee^{B-x} ,\quad &\text{if } j\leq  J,
\\
  K\, \ee^{ \alpha_n-\alpha_{n- j } } \ee^{-j^{\varrho}},\quad &\text{if } j\in (J ,\frac{2}{3}n],
  \\
  \ee^{ \alpha_{n } } \ee^{- j^{ \varrho }},\quad &\text{if } j\in (\frac{2}{3}n,t] .
\end{array}\right.
\end{eqnarray*}

By summing these inequalities,  for $n$ large enough we get that  
\begin{eqnarray*}
\mathbb Q_{\eqref{goodwn}} 
&\le&   K' (J \,  \ee^{B-x}  + \ee^{- J^{\varrho/2}}) \sum_{t=n-T}^n    \q\Big(  F_n , \, t=\tau_{\zeta_n}({w_n})  ,\,  w_n \text{ good}     \Big)  
	\\ &\le &  K' (J \,  \ee^{B-x}  + \ee^{- J^{\varrho/2}})\sum_{t=n-T }^n    \P \left( S_n\leq  \alpha_n -x ,\,  \min_{\tau_{\zeta_n} \le i \le n}S_i \geq  \alpha_n -x-L,\, t=\tau_{\zeta_n},\, \tau_{\zeta_n}^{(2)}>n   \right)  
\\
&\leq & K'' (J \,  \ee^{B-x}  + \ee^{- J^{\varrho/2}}) (1+L) \; T\,   \ee^{-\alpha_n},  
\end{eqnarray*}
where for the second inequality we have used \eqref{manytoone} for $F_n \cap\{ t=\tau_{\zeta_n}({w_n})  \}$ and for the last inequality, we have applied  (\ref{pasprov}) to      $y= \alpha_n-x-L$ and $a= 1, ..., \lceil L\rceil$).   Finally we choose  $J=J(L, T, K'')$ large  enough and $x\ge x_1(B, J)$ so  that $K'' (J \,  \ee^{B-x}  + \ee^{- J^{\varrho/2}}) (1+L) \; T\, \le \varepsilon$.   Then   $\q_{\eqref{goodwn}}  \le \varepsilon \ee^{-x}$ and \eqref{eqclust} follows.  This proves Proposition \ref{clust}.   
\hfill$\Box$

\medskip

We end this section by the proof of Lemma \ref{L:App1}.

\medskip
\noindent{\it Proof of Lemma \ref{L:App1}.}   
%The event  ${\mathbb G}_{n, J, B}^c$means that there exists $k\in [1,\tau_{\zeta_n}^{(w_n)}]$ such that (\ref{defGood}) is not satisfied. 

Firstly we will prove that with overwhelm probability the trajectory of $(V(w_i))_{i\geq 0}$ contains only one big jump and never drops too low.    Recall the notation $F_n$ defined in \eqref{def-fn}. 
Write for brevity  $$
	y :=  \alpha_n - x- L.$$

We shall use several times the fact that under $\q$, $(V(w_j), j\ge 0)$ has the same law as $(S_j, j\ge0)$ under $\P$. Then by   (\ref{pasprov2bis}) with $a=0, 1, ..., \lceil L \rceil$, we get that for some constant $K_L>0$ depending on~$L$,  
\begin{eqnarray}
 \q\left(F_n , \tau_{\zeta_n}^{(2, w_n)}\leq n  \right)
 \nonumber &&=  \sum_{i=n-T}^n \P \left( \underline{S}_{ [\tau_{\zeta_n},n]}  \geq  y ,\, \tau_{\zeta_n} =i,\, S_n \in [y, y +L],\, \tau_{\zeta_n}^{(2)}\leq n \right) 
\\
\label{elag1} &&\leq  K_L\, \sum_{i=n-T}^n  n^{-2\alpha} i^{-1/\alpha} \ell_3(n)   \leq \varepsilon   \, \ee^{-\alpha_n} ,
\end{eqnarray}

\noindent
for all large  $n$.  We claim  that there exists some positive constant $c_4=c_4(L,T)$ such that for all $n$ large,  \begin{equation} \label{elag2}
\q\Big(F_n,\, \min_{1\le j <\tau_{\zeta_n}^{(w_n)}} V(w_j) \leq -c_4, \, \tau_{\zeta_n}^{(2, w_n)}> n  \Big) \le  O(  \varepsilon )  \, \ee^{-\alpha_n} . 
\end{equation}

Let us denote by $\q_{\eqref{elag2}}$ the probability term in \eqref{elag2}.   Denote by $j$ be the first time such that $V(w_j) \le - c_4$; then 
by using  the Markov property at $j$,  we get that
\begin{eqnarray*}
\q_{\eqref{elag2}}&  =  &\sum_{i=n-T}^n\sum_{j=1}^{i-1}\E\Big[ 1_{\{  \underline{S}_{j-1} > -c_4,\, S_j\leq -c_4,\, \underset{k\leq j}{\min}X_k\geq -\zeta_n\}}  \times
\\
&&\qquad\qquad\qquad\qquad\qquad \P_{S_j }\Big(S_{n-j}\leq y+L  ,\, \min_{\tau_{\zeta_n} \le i \le n-j}S_i \geq y ,\, \tau_{\zeta_n}=i-j ,\, \tau_{\zeta_n}^{(2)}> n-j \Big) \Big]
\end{eqnarray*}
If $j\geq \frac{n}{2}$ by Lemma \ref{L:1} (with $x= \mm j +c_4$ and $y= \zeta_n$ there) we get that  $$\P\big( S_j\leq -c_4,\, \underset{k\leq j}{\min}X_k\geq -\zeta_n  \big)\leq K\ee^{- (\log n)^3/K},$$

\noindent whereas  for  $j\leq \frac{n}{2}$, since   $y- S_j \le y+ c_4+ \zeta_n \leq {\mm\over 2} (n-j)  $, by using $\lceil L\rceil $ times (\ref{pasprov}) (with      $a   \in [0,L]$  being integer),  we deduce that   
for any $i\in [n-T,n]$ and on $\{S_j\leq -c_4,\, \underset{k\leq j}{\min}X_k\geq -\zeta_n\}$, 
\begin{equation}
  \P_{S_j }\left(S_{n-j}\leq y+L  ,\, \underline{S}_{[\tau_{\zeta_n},n-j]}\geq y ,\, \tau_{\zeta_n}=i-j ,\, \tau_{\zeta_n}^{(2)}> n-j \right) \leq   K' \, (1+L)  \, \ee^{-\alpha_n} 
 \end{equation}
(we used the fact that  $\ee^{-\alpha_{n-j}}=O(\ee^{-\alpha_{n}})$ when $j\le \frac{n}{2}$). It follows   that 
\begin{eqnarray*}
\q_{\eqref{elag2}}& \le&  \sum_{i=n-T}^n\sum_{j=\frac{n}{2}}^{i-1} K\ee^{- (\log n)^3/K} + K'\, (1+L) \, \ee^{-\alpha_n} \, \sum_{i=n-T}^n\sum_{j=1}^{\frac{n}{2}}  \P\Big( \underline{S}_{j-1}\geq - c_4,\, S_j<-c_4  \Big)  \\
  &&\leq \varepsilon  \, \ee^{-\alpha_n}   +K'\, (1+L) \, T\,    \ee^{-\alpha_n} \,      \sum_{j=1}^\frac{n}{2} \P\left(\underline{S}_{j-1}\geq -c_4,\, S_j<- c_4 \right)
\\
  && \leq \varepsilon   \, \ee^{-\alpha_n} + K'\, (1+L) \, T\,  \, \ee^{-\alpha_n} \,    \P\left(\min_{k\geq 0}S_k <  -c_4\right) \\
  && \leq 2\varepsilon    \, \ee^{-\alpha_n} ,
\end{eqnarray*}

\noindent 
by  choosing $c_4=c_4(L, T)$   large enough to get the last inequality.  Then  \eqref{elag2}  follows.

 \medskip
 
By  combining (\ref{elag1}) and (\ref{elag2}),  to get Lemma \ref{L:App1} it is enough to prove the following assertion: {\it for any $L,\, T,\, \varepsilon>0$ there exist  $B>0$ and $J$  such that for any $n\ge n_0(J, B, L, T, \varepsilon)$, $ x\geq 0$,
\begin{equation} \label{last6} 
\q\Big(F_n,   \,   \min_{1\le j < \tau_{\zeta_n}^{(w_n)}} V(w_j)  \geq -c_4,\,  \tau_{\zeta_n}^{(2, w_n)}>n,\,   \mbox{ $w_n$ not good} \Big)\leq  O(\varepsilon) \, \ee^{-\alpha_n} . 
\end{equation}  }

Recall that  $0<\varrho<\min(\frac{\alpha-1}{2},\frac{1}{12})$.  Before establishing \eqref{last6} we prove the following claim:

\begin{Claim} \label{Claim52}

(i)    There is   a sequence of positive real numbers $(\varepsilon_j)$ such that $\lim_{j\to \infty} \varepsilon_j=0$ and for any integer $j$ and   $z\in \r , \, y\geq 0$,
\begin{equation} \label{App2}
\sum_{p= j}^\infty \P\left(S_p-  p^{\varrho}\leq z,\, \underline{S}_p\geq -y \right) \leq  \frac{10}{ \mm}\, \big(z-\frac{\mm j }{10}\big)^+  + \varepsilon_j (1+y+z^+).\end{equation} 

(ii)  There  exits some positive constant $K_{L, T} >0$ such that for all large $n$  and  $k\in [1,n-T)$,  \begin{equation}
\sup_{z \le -  \frac{mn}{5}} \P \left(S_{n-k}\leq   z,\,  \min_{\tau_{\zeta_n}\le i \le n-k} S_i \geq z -L,\,  \tau_{\zeta_n}  \in [n-k-T,n-k],\, \tau_{\zeta_n}^{(2)}> n-k  \right) \le   K_{L, T} \, \ee^{-\alpha_n} .  \label{newplus} 
\end{equation}
\end{Claim}

\noindent{\it Proof of Claim \ref{Claim52}.}  

(i) Observe that
\begin{eqnarray*}
\sum_{p\geq j}\P\left(S_p-  p^{\varrho}\leq z ,\, \underline{S}_p\geq -y\right)&\leq&  \Big (\frac{10}{m}z- j\Big)^+   + \sum_{p\geq  \max(\frac{10 z}{m},  j)}\P\left(S_p \in [-y, z+p^{\varrho}] \right).
\end{eqnarray*}

Then observe that $z+p^{\varrho} \le  {\mm\over 2} p$ for all $p \ge  {10 z \over m}$ and $p\ge j_0$ if $j_0$ is large enough.  By applying Lemma~\ref{3urd}, we get 
 \begin{eqnarray*}
 \sum_{p\geq j}\P\left(S_p-  p^{\varrho}\leq z ,\, \underline{S}_p\geq -y\right)&\leq& \frac{10}{ \mm}\, \big(z-\frac{\mm j }{10}\big)^+    + K\, \sum_{p\geq  \max(\frac{10 z}{m},  j)} \frac{l(p)(y+ z^++p^{\varrho})}{p^\alpha}
 \\
&  \leq & \frac{10}{ \mm}\, \big(z-\frac{\mm j }{10}\big)^+  + \varepsilon_j (1+y+z^+) ,
 \end{eqnarray*}
 with $\varepsilon=O(j^{(1-\alpha)/2})$, proving \eqref{App2}.
 
\medskip

 (ii) Denote by $\P_{\eqref{newplus}}$ the probability term in \eqref{newplus}.  Then $$ \P_{\eqref{newplus}} = \sum_{j=n-k-T}^{n-k} \P \left(S_{n-k}\leq   z,\,  \min_{j \le i \le n-k} S_i \geq z -L,\,  \tau_{\zeta_n} =j,\, \tau_{\zeta_n}^{(2)}> n-k  \right)=: \sum_{j=n-k-T}^{n-k}  \P_{\eqref{newplus}} (j) .$$
 
 Notice that  $z-L \le S_{n-k}\leq   z$. Therefore if   $S':= S_{n-k} - X_j \ge -{\mm\over10} n$ then $X_j \le z +  {\mm\over10} n \le -{\mm\over10} n$.  Moreover $ z-L - S' \le X_j \le z - S'$. By   the independence of $X_j$ and $S'$, we get that   $$ \P \left(z-L \le S_{n-k}\leq   z,\,  \tau_{\zeta_n} =j,\,S'  \ge -{\mm\over10} n  \right)   \le    \sup_{ b \le   -{\mm\over10} n} \P\left( X_j \in [b-L, b]\right) \le K_L\, \ee^{-\alpha_n},$$ 
 
 \noindent by using the density of $X_j$ given by  \eqref{Hyp2}.  On the other hand, if $S':= S_{n-k} - X_j <  -{\mm\over10} n$ then we can apply Lemma \ref{L:1} to see that $$\P \left(S'<  -{\mm\over10} n,\,  \tau_{\zeta_n} =j,\, \tau_{\zeta_n}^{(2)}> n-k  \right) \le \P \left( S_{n-k-1} <  -{\mm\over10} n, \tau_{\zeta_n} > n-k-1\right)  \le K\, \ee^{- \mm n / (10 \zeta_n)}.  $$

 \noindent Therefore $ \P_{\eqref{newplus}} (j)  \le K_L\, \ee^{-\alpha_n}+ K\, \ee^{- \mm n / (10 \zeta_n)}$ and \eqref{newplus} follows if we take $K_{L, T}= 2 (1+T) K_L.$  This completes the proof of Claim \ref{Claim52}. \qed

\medskip
Let us go back to the proof of \eqref{last6}.    Define  for any $k \geq 1$, $\xi(w_k):=  \sum_{v\in {\mathbb B}(w_k)} \ee^{- \Delta V(v) }$. Then,    \begin{equation}
\label{consdefxi}\sum_{v\in {\mathbb B}(w_k)}\ee^{-(V(v)+x)}= \ee^{-V(w_{k-1})-x}\xi(w_k). 
\end{equation}

\noindent Notice  that the sequence $\{\xi(w_k),  \Delta V(w_k) \}_{k\ge1}$ are i.i.d. under $\q$.  Define $\xi= \xi(w_1)$.  

Let $n$ be large enough so  that $n-T> J$. On $\{ \tau_{\zeta_n}^{(2, w_n)}>n\} \cap  \{ \mbox{ $w_n$ not good}\} $,    either there is some $1\le k \le  J$ such that $\xi(w_k) > \ee^{ B + V(w_{k-1})}$ or some $  J < k  <  \tau_{\zeta_n}^{(w_n)}$ such that $\xi(w_k) > \ee^{ V(w_{k-1})+x - k^{\varrho}}$.   We   discuss separately  these two  cases:

\medskip
{\noindent\bf  The first case: choice of $J$ to control  $J  <  k  < \tau_{\zeta_n}^{(w_n)}   $ :}  

Notice that $\tau_{\zeta_n}^{(w_n)}\in [n-T, n]$. For any   $J < k < n-T$, 
 we apply  the Markov property at $k$ to arrive at  
\begin{eqnarray}
 \q_{\eqref{hehe1}} (k) &:=& \q\Big(k  < \tau_{\zeta_n}^{(w_n)}, \,  \xi(w_k) > \ee^{x+V(w_{k-1})-k^{\varrho}},\,  \min_{1\le j \le  k}V(w_j) \geq -c_4,\,  F_n , \, \tau_{\zeta_n}^{(2, w_n)}> n  \Big) \nonumber
\\
&  =&\mathbb E_\q\Big[ 1_{\{ \xi(w_k) > \ee^{x+ {V(w_{k-1})-k^{\varrho}}},\, \min_{1\le j \le  k}V(w_j) \geq - c_4,\, \underset{j\leq k}{\min} \Delta V(w_j) \geq -\zeta_n  \}}    \, g_{n-k}(V(w_k))     \Big], \label{hehe1}
\end{eqnarray}

\noindent with $$ g_{n-k}(b):= \P \Big( S_{n-k}\leq y -b+L,\,  \min_{\tau_{\zeta_n} \le j \le n-k}S_j \geq y -b ,\,  \tau_{\zeta_n} \in [n-k-T,n-k],\, \tau_{\zeta_n}^{(2)}> n-k \Big), \, \,  z \in \r. $$

When $k\leq \frac{n}{2}$,  we can apply (\ref{pasprov-sum})  to get that  for $b:= V(w_k) \ge - c_4$, $$ g_{n-k} (b) \le K'_L\, \ee^{-\alpha_n},$$
with some positive constant $K'_L$ depending on $L$ [in fact $K'_L =O(L^2)$].  

\medskip

 For $k\geq \frac{n}{2}$,     if  $ b:= V(w_k) >  {\mm n\over 4}$, then  $  y+L - b \le  {-\mm n \over5}$ and $g_{n-k}(b) \le  K_{L, T} \ee^{-\alpha_n}  $ by \eqref{newplus}.  Consequently, we get that for any $J < k <n -T$ and $x\ge0$, $$  \q_{\eqref{hehe1}} (k) \le K''_{L, T} \ee^{-\alpha_n} \, \q \Big( \xi(w_k) > \ee^{ {V(w_{k-1})-k^{\varrho}}}\Big) + 1_{\{ k \ge {n\over2}\}}  \q \Big( V(w_k) \le  {\mm n\over 4},   \underset{j\leq k}{\min} \, \Delta V(w_j) \geq -\zeta_n\Big). $$
Moreover, $ \q \left( V(w_k) \le  {\mm n\over 4},   \underset{j\leq k}{\min} \, \Delta V(w_j) \geq -\zeta_n\right)= \P\left( S_k \le  {\mm n\over 4},   \underset{j\leq k}{\min} \,   X_j  \geq -\zeta_n\right) \le K \, \ee^{- \mm n /(4\zeta_n)}$ by Lemma \ref{L:1}. Since under $\q$, $\xi(w_k)$ is independent of $V(w_{k-1})$ which is distributed as $S_{k-1}$ under $\P$, and, moreover, $\xi(w_k)$ has the same law as $\xi$,  it follows from  (\ref{App2})  that \begin{eqnarray*}
\sum_{J< k \le n-T }  \q_{\eqref{hehe1}} (k) & \le & K''_{T, L}\, \ee^{-\alpha_n} \,  \mathbb E_\q\Big[ \sum_{J< k \le n-T }  \P( \log \xi(w_k) \geq  S_{k-1}-k^{\varrho},\, \underline{S}_{k-1}\geq - c_4 \})  \Big] + K \,n \,  \ee^{- \mm n /(4\zeta_n)}  \\
	&\le& K_{T, L}\, \ee^{-\alpha_n} \, \Big( \mathbb {E_Q}\Big[ (\log \xi-\frac{\mm J }{10})^+\Big] +  \varepsilon_J \,\mathbb {E_Q} \big[ 1+  c_4+( \log \xi)^+\big]\Big)+ K \,n \,  \ee^{- \mm n /(4\zeta_n)} ,
\end{eqnarray*}

\noindent with $\varepsilon_J \to 0$ as $J \to \infty$.  By  (\ref{Hyp3}), $\mathbb{E_Q}[(\log \xi)^+ ]\leq \mathbb E \Big[\sum_{|u|=1} \ee^{-V(u)} \big(\log [\sum_{|u|=1} \ee^{-V(u)} ]\big)^+\, \Big]<\infty$, thus we choose and then fix $J=J(\varepsilon, T, L)$ large enough so that $  \mathbb {E_Q}\left[ (\log \xi-\frac{\mm J }{10})^+\right] +  \varepsilon_J \,\mathbb {E_Q} \big[ 1+  c_4+( \log \xi)^+\big]\Big)  \le {\varepsilon\over K_{L, T}} $. 
Then for all large $n$, we get that 
 \begin{equation}\label{first98} \sum_{J< k \le n-T }  \q_{\eqref{hehe1}} (k) \le 2 \varepsilon \, \ee^{-\alpha_n}.\end{equation}

{\noindent \bf  The second (and last) case:  Choice of $B$ to control   $1\le k \leq  J$: } 

 Under $\q$ and conditionally on $\{V(w_k)=z\}$, the process $\{V(w_{i+k}), 0\le i \le n-k\}$ is distributed as $\{ S_i, 0\le i \le n-k\}$ under $\P_z$.  It follows from   the Markov property  at $k$ that 
\begin{eqnarray*}
&&\q\Big( \xi(w_k) \geq \ee^{B+V(w_{k-1})}  ,\,  F_n  ,\, \min_{1\le j \le  k}V(w_j) \geq - c_4,\,   \tau_{\zeta_n}^{(2, w_n)}> n \Big)
\\
&&=\mathbb E_\q\Big[ 1_{\{ \xi(w_k) \geq \ee^{B+V(w_{k-1})},\,   \min_{1\le j \le  k}V(w_j) \geq - c_4 \}} \times
\\
&&\qquad\qquad\qquad   \P_{V(w_k) }\left( S_{n-k}\leq y+L,\,  \min_{\tau_{\zeta_n}-k \le j \le n-k}S_j \geq y ,\,  \tau_{\zeta_n} \in [n-k-T,n-k],\, \tau_{\zeta_n}^{(2)}> n-k  \right)\Big]
\\
&&\leq  K_L    \, {\q}\left( \xi(w_k) \geq \ee^{B+V(w_{k-1})},\, \min_{1\le j \le  k}V(w_j) \geq - c_4 \right) \, \ee^{-\alpha_n}   ,
\end{eqnarray*}
where $K_L>0$ denotes some constant depending on $L$ and  we have applied    (\ref{pasprov-sum}) to get  the last inequality  [remark that  $y-V(w_k) \le y+c_4\le {\mm\over2}(n-k)$]. Furthermore, 
\begin{eqnarray*}
  \mathbb E_\q\Big[  \sum_{k=1}^{J}  1_{\{ \log \xi(w_k) \geq B+V(w_{k-1})  ,\,\min_{1\le j \le  k}V(w_j) \geq -  c_4\}}  \Big] 
	 & \leq   &\sum_{k=1}^{J} \q\left( V(w_{k-1})  \leq -\frac{B}{2} \right)  +   \sum_{k=1}^{J}\frac{2}{B}\mathbb E_\q\big[ (\log  (\xi(w_k)))^+\big]
	\\ &  =  & \sum_{k=1}^{J} \P\left( S_{k-1}\leq -\frac{B}{2} \right)  +  \frac{2 J}{B} \mathbb E_\q\big[ (\log \xi)^+\big]
	\\ &\leq & {\varepsilon\over K_L}, 	\end{eqnarray*}

\noindent by choosing $B=B(J, L, T, \varepsilon)$ large enough. Finally we have 
 \begin{equation}\label{first5} \q\Big( \exists k\in [1, J]: \xi(w_k) \geq \ee^{B+V(w_{k-1})}  ,\,\min_{1\le j \le  k}V(w_j) \geq - c_4,\,   F_n  ,\, \tau_{\zeta_n}^{(2, w_n)}> n \Big) \le \varepsilon \, \ee^{-\alpha_n}.\end{equation}

 By combining  \eqref{first98} and \eqref{first5}, we get \eqref{last6} and therefore complete the proof of   Lemma \ref{L:App1}.\hfill$\Box$

%{\bf Case 3: $k  = \tau_{\zeta_n}^{(w_n)} $.}    It remains to estimate \begin{eqnarray*} \nonumber &&\q\Big( \sum_{v\in {\mathbb B}(w_{\tau_{\zeta_n}({w_n})})} \ee^{-(V(v)+x)} \geq \ee^{-k^\frac{1}{12}},\,   V(w_n)\leq \alpha_n-x ,\, \underline{V}_{[\frac{n}{2},n]}(w_n)\geq \alpha_n-x -L,\, \tau_{\zeta_n}({w_n})\in [n-T,n] \Big) \\ \nonumber &&\leq \sum_{t=n-T}^n \q\left(  \xi(w_t) \geq \ee^{V(w_{t-1})-k^{\frac{1}{12}}},\,   V(w_n)\leq \alpha_n-x ,\, \underline{V}_{[\frac{n}{2},n]}(w_n)\geq \alpha_n-x -L,\, \tau_{\zeta_n}({w_n})=t \right) \\ &&\leq \varepsilon \mm^{-(\alpha+1)} \, \ee^{-\alpha_n}  \end{eqnarray*} According to the hypothesis (\ref{decorr})( A PRECISER ). 

 %-Let $\delta>0$. By a trivial decomposition we deduce that for any $j\geq 1$, 
%\begin{eqnarray}
%\nonumber \P\left( S_j-mj \leq -\delta j\right) &= &\P(S_n-mj<-\delta j,\, \underset{1\leq i\leq n}{\min} X_i\geq -\frac{j}{(\log j)^3})+ \P(S_j-mj<-\delta j,\, \underset{1\leq i\leq n}{\min} X_i\leq -\frac{j}{(\log j)^3})
%\\
%\label{2urd}&\leq & C\ee^{-\delta(\log j)^3}+ c\frac{j}{j^{\alpha-\varepsilon}}\leq \frac{C}{j^{\alpha -1-\varepsilon}}.
%\end{eqnarray}

%-For any $n$ large enough
%\begin{eqnarray}
%\nonumber \sup_{y\in \r}\P\left( S_n-y \leq [0,1],\, \tau_{\zeta_n}(S)=n-i\right) &\leq &\P\left(X_{n-i}\leq -\zeta_n,\, X_{n-i}+ S_{n-i-1}+ S_{n}-S_{n-i-1}\in [0,1]\right)
%\\
%\label{biss} &\leq&  \sup_{x\leq -\zeta_n}\P\left( X_{n-i}\in [x,x+1]\right)\leq \frac{l(\zeta_n)}{(\zeta_n)^{\alpha+1}}
%\end{eqnarray}

%\bibliographystyle{plain}
%\bibliography{bibli}

\end{document}